\documentclass [10pt,reqno]{amsart}
\usepackage{ucs}
\usepackage[latin1]{inputenc} 
\usepackage{amsmath, amssymb, amscd}
\usepackage{pb-diagram}
\usepackage[latin1]{inputenc}
\usepackage{amsmath}
\usepackage{amssymb}
\usepackage{hyperref}

\newtheorem{theorem}{Theorem}[section]
\newtheorem{definition}[theorem]{Definition}
\newtheorem{lemma}[theorem]{Lemma}
\newtheorem{prop}[theorem]{Proposition}
\newtheorem{proposition}[theorem]{Proposition}
\newtheorem{corollary}[theorem]{Corollary}

\newtheorem{condition}[theorem]{Condition}

\newtheorem*{observation}{Observation}
\newtheorem{remark}[theorem]{Remark}

\newcommand{\abs}[1]{\lvert#1\rvert}
\newcommand{\norm}[1]{\| #1 \|}

\renewcommand{\epsilon}{\varepsilon}
\newcommand{\dirac}{\diagup \hskip-9pt D}

\DeclareMathAlphabet{\mathpzc}{OT1}{pzc}{m}{it}

\usepackage{mathrsfs}

\newcommand{\N}{\mathbb{N}}
\newcommand{\Z}{\mathbb{Z}}
\newcommand{\C}{\mathbb{C}}

\newcommand{\R}{\mathbb{R}}

\renewcommand{\qed}{$\hfill \square$ \medskip \\}
\renewcommand{\phi}{\varphi}

\newcommand{\ad}{\text{ad}}

\newcommand{\SO}{SO}

\newcommand{\tensor}{\otimes}

\newcommand{\Hom}{\text{Hom}}
\newcommand{\Aut}{\text{Aut}}
\renewcommand{\det}{\text{det}}
\newcommand{\coker}{\text{coker}}
\newcommand{\im}{\text{im}}
\newcommand{\id}{\text{id}}
\newcommand{\Hol}{\text{Hol}}
\newcommand{\ind}{\text{ind}}

\newcommand{\A}{\mathscr{A}}

\newcommand{\su}{\mathfrak{su}}
\renewcommand{\sl}{\mathfrak{sl}}
\newcommand{\gl}{\mathfrak{gl}}

\newcommand{\G}{\mathscr{G}}

\newcommand{\data}{\mathfrak{s},E}

\newcommand{\conf}{\mathscr{C}}
\newcommand{\bonf}{\mathscr{B}}
\newcommand{\V}{\mathscr{V}}
\newcommand{\W}{\mathscr{W}}
\newcommand{\s}{\mathfrak{s}}

\newcommand{\bfomega}{\boldsymbol{\omega}}
\newcommand{\bfalpha}{\boldsymbol{\alpha}}
\newcommand{\bfbeta}{\boldsymbol{\beta}}
\newcommand{\bfnu}{\boldsymbol{\nu}}

\newcommand{\stab}{\Gamma \times_{\Z/2} S^1}

\begin{document}
\thispagestyle{empty}
\title{A vanishing result for a Casson-type instanton invariant} 
\author[Raphael Zentner]{Raphael Zentner} 
\address{Mathematisches Institut, Universit\"at zu K\"oln, Weyertal 86--90, 50931 K\"oln \\
Germany}
\email{rzentner@math.uni-koeln.de}
\maketitle
\begin{abstract}
Casson-type invariants emerging from Donaldson theory over certain negative definite 4-manifolds were recently suggested by Andrei Teleman. These are defined by an algebraic count of points in a zero-dimensional moduli space of flat instantons. Motivated by the cobordism program of proving Witten's conjecture, we use a moduli space of $PU(2)$ Seiberg-Witten monopoles to exhibit an oriented one-dimensional cobordism of the instanton moduli space to the empty space. The Casson-type invariant must therefore vanish. 
\end{abstract}
\subsection*{MSC 2010:} 57R57, 57N13, (53C27, 58D27, 58J37)

\section*{Introduction}
In a recent work of Andrei Teleman \cite{T} low energy instanton moduli spaces defined over smooth, closed, oriented, negative definite four-manifolds $X$ appeared. All these moduli spaces are compact and do not contain reducibles. Among these are certain `Casson-type' moduli spaces defined for manifolds  with $b_2(X) \equiv 0 \ (\text{\em mod} \ 4)$. If in addition $b_1(X) = 1$ the expected dimension of the moduli space is zero, and a `Casson-type' invariant can be defined by an algebraic count of the elements in the moduli space, provided the latter is regular. This type of invariant has been suggested by Teleman, although an explicit definition of the invariant is missing in \cite{T}. 

Andrew Lobb and the author have shown that non-emptiness of the Casson-type moduli space also gives an
obstruction to certain connected sum decompositions of the four-manifold, ruling out natural appearing
candidates for a non-empty moduli space \cite{LZ}. We also have studied a method to obtain manifolds with non-empty
Casson-type moduli space by surgery on  regular neighbourhoods of two-knots in $4 \overline{\mathbb{CP}^2}$
with trivial normal bundle. We are confident that this method will eventually yield non-empty moduli spaces \cite{LZ}, at least for manifolds with boundary.

Similar moduli spaces and Casson-type invariants have been defined by Ruberman and Saveliev for $\Z[\Z]$-homology Hopf surfaces \cite{RS}, going
back to work of Furuta and Ohta \cite{FO}, and for $\Z[\Z]$-homology 4-tori \cite{RS2}. The first of these is also related in a non-trivial way to Seiberg-Witten theory \cite{MRS} by work of Mrowka, Ruberman and Saveliev. In particular, these invariants are in general non-vanishing. Non-vanishing results were also expected for the Casson-type invariants considered here. 
However, the author was informed in June 2008 by Andrei Teleman that in the perspective of the `cobordism program' for proving Witten's conjecture he conjectured the invariant to vanish. He expected that this could be proved using the formalism of `virtual fundamental classes'. The author, however, saw the possibility of proving this with a genuine cobordism provided one uses suitable (holonomy) perturbations. Even though the idea is heuristically easily described making it precise requires a lot of technical work. In particular, the machinery developed by Feehan and Leness \cite{FL2,FL3,FL,FL_general}, and the results by Teleman \cite{T3} and Okonek-Teleman \cite{OT} do not apply to the case required for our situation, at least not directly.
\\

In this article we will give an explicit definition of the Casson-type instanton invariant suggested by Teleman. We will actually show that his idea can be extended to a more general situation, where the requirement on the second Betti number is weakened to $b_2(X)\geq 4$ (thanks to a remark by Kim Fr\o yshov made to the author). We will generalise the holonomy perturbations as used by Peter Kronheimer in \cite{K} for higher rank instanton moduli spaces in order to obtain generic regularity for the monopole moduli spaces we consider. Natural orientations for the moduli space, obtained by the determinant line bundle of the (Fredholm-) deformation operator parametrised by the configuration space, then yields the Casson-type invariant as algebraic count of points in the perturbed moduli space. That this signed count is independent of the chosen perturbation then follows from the fact that the parametrised moduli space doesn't contain reducibles, and the fact that the orientation was chosen in the given natural way from the ambient configuration space containing the moduli space.
\\


The idea of the cobordism program is to use $PU(2)$ monopole moduli spaces containing both the instanton moduli space defining the Donaldson-invariant, and certain Seiberg-Witten moduli space as fixed point subspaces of a natural circle-action. The $S^1$-quotient of this moduli space then yields a cobordism between a bundle over the instanton moduli space and bundles over the Seiberg-Witten moduli spaces, at least heuristically, provided that the index of the coupled Dirac operator is positive. There are many very difficult technical problems involved, mainly because the moduli space needs not to be compact and parts of the instanton moduli space and some Seiberg-Witten moduli spaces may lie in lower strata of the Uhlenbeck-compactification. 

In our situation the Casson-type instanton moduli space is compact and zero-dimensional. We will show that a compact two-dimensional $PU(2)$ monopole moduli space can be found, and that the $S^1$ quotient yields an oriented cobordism between the Casson-type moduli space and the empty space. This implies that the algebraic count is zero, so the Casson-type invariant vanishes. 
\\

Technically our approach to transversality is essentially different from Feehan and Leness's who have used holonomy perturbations as well in earlier versions of their work \cite{FL2}. In their situation the holonomy pertubation term of the curvature equation becomes zero if the spinor vanishes. This is fine if one still can make this locus regular, for instance, in their situation by perturbing the Riemannian metric. However, this is not possible in our situation as we will have an instanton moduli space of flat connections. Secondly, their holonomy perturbations are supported in a finite number of balls on the manifold, and they require elaborate arguments that imply that local reducibility for the connection component of the monopoles imply global reducibility. This continuation theorem seems to fail for instantons, see the Remark 5.24 in \cite{FL2}, and so we cannot apply their technique to get generic regularity of the instanton moduli space. Apparently, their `indirect' way to kill the cokernels of the derivative of the parametrised monopole map is necessary to allow the development of glueing theory in their program of proving Witten's conjecture.

In our situation we cannot allow the holonomy perturbation term of the curvature equation to vanish if the spinor component vanishes. On the other hand, we need to add somewhat artificial cut-off factors in the perturbed monopole equations in order to obtain an a-priori bound on the spinor component. As we do not have control of this cut-off factor for a sequence of monopoles that `goes to infinity' we use a second cut-off factor that vanishes for a sequence tending to infinity. As a result, the monopoles (and instantons) appearing in the lower strata appear all without holonomy perturbations. There we use classical perturbations to obtain regularity, following Teleman \cite{T}. Similarly, classical perturbations were also studied by Feehan in \cite{F}.  
\\

From the point of view of investigating invariants of smooth four-manifold our vanishing result is negative. However, it appears that this result has relevance to Teleman's classification program on complex surfaces of class VII \cite{T2,T4}, and in this perspective, it is rather a positive result. We discuss this briefly in Section \ref{Teleman's program} below. \\

\section*{Acknowledgements}
The author would like to thank Kim Fr\o yshov, Daniel Ruberman, and Andrei Teleman for helpful discussions on these Casson-type instantons and related gauge theory. He is also very grateful for an anonymous referee's detailed mathematical and linguistic comments, and for a second referee's further comments.

\section{Casson-invariant for negative definite four-manifolds}
In this section we define Teleman's Casson-type instanton invariant via gauge theory of anti-self-dual connections on appropriate bundles. We shall first define the configuration space and fix our notation conventions, then define the moduli space and discuss some of its properties. Then a suitable space of holonomy perturbations is introduced in order to get a regular moduli space. After introducing preferred orientations for the zero-dimensional moduli space we define the invariant. The book of Donaldson and Kronheimer \cite{DK} can be seen as a general reference here.

\subsection{The configuration space}\label{elementary definition configuration space} 
Let $X$ be a smooth closed oriented Riemannian four-manifold. Let $E \to X$ be a Hermitian rank-2 bundle. We suppose a smooth connection $a$ on the determinant line $\det(E)$ is fixed, and define the space $\mathscr{A}_a(E)$ to be the affine space of unitary connections in $E$ which induce the fixed connection $a$ in the determinant line bundle,
 and which are of Sobolev class $L^2_l$, for some fixed $l \geq 3$.

The choice of $l$ to be greater or equal than $3$  ensures that there is a Sobolev embedding $L^2_l \hookrightarrow C^0$ which we will use further down in the Transversality Theorem \ref{transversality monopoles}. We define the `gauge group' $\mathscr{G}$ to be the group of (unitary) automorphisms of $E$ with determinant one and of Sobolev class $L^2_{l+1}$. It acts by the formula $(u,\nabla_A) \mapsto u \circ \nabla_A \circ u^{-1}$ if we denote by $\nabla_A$ the covariant derivative of the linear connection $A$ and $u \in \G$.  The quotient $\mathscr{B} = \mathscr{A}_a(E) / \G$ is called the configuration space. We shall denote by $\Gamma_A \subseteq \G$ the stabiliser of the connection $A$ under the gauge group $\G$. It is the centraliser of the holonomy group associated to $A$. The centre $Z = \Z/2$ of $SU(2)$, seen as constant gauge transformations on $X$, is in the stabiliser $\Gamma_A$ for any connection. A connection $A$ is called irreducible if $\Gamma_A = Z$, otherwise reducible. A connection $A$ on $E$ is reducible if and only if there is a proper $A$-invariant subbundle of $E$. 

The following is an equivalent viewpoint of the space $\mathscr{A}_a(E)$ and the configuration space $\mathscr{B}$: Denote by $P \to X$ a principal bundle with structure group $U(2)$ and by $\ad(P) \to X$ the $PU(2)$ bundle associated to $P$ by the adjoint representation. Let $\mathscr{A}(\ad(P))$ denote the space of $PU(2)$ connections in $\ad(P)$ of class $L^2_l$. The group of automorphisms of $P$ of determinant one and Sobolev class $L^2_{l+1}$ acts on $\mathscr{A}(\ad(P))$. The space $\mathscr{A}(\ad(P))$ is naturally isomorphic to the space $\mathscr{A}_a(E)$ above if the bundle $P$ is taken as the unitary frame bundle $P_E$ of $E$, and this is equivariant with respect to the mentioned automorphism groups. In this case the bundle $\ad(P)$ is naturally identified with the bundle $\su(E)$ of traceless skew-adjoint endomorphisms of $E$. Unless we want to make explicit our viewpoint we shall simply write $\mathscr{A}$ for the spaces of connections we have in mind, or  we shall write $\mathscr{A}_l$ for this space if we want to make the dependence on the Sobolev index explicit. 

Principal $U(2)$ bundles or Hermitian vector bundles of rank 2 on a manifold of dimension 4 are classified by their first and second Chern-classes. However, it is a common convention in gauge theory, notably in our main references \cite{DK}, \cite{FL}, \cite{K} to encode this information in the determinant line bundle $w \to X$ (or its first Chern class) of $E$ respectively $P$, and in the `instanton number' 
\begin{equation}\label{instanton number}
 \kappa = -\frac{1}{4} \langle p_{1}(\su(E)),[X] \rangle = \langle c_{2}(E) - \frac{1}{4} c_1(E)^2, [X] \rangle \ .
\end{equation}
In notation we will stick to the common convention and write $\mathscr{B}^w_\kappa$ for the configuration space $\mathscr{B}$ if we want to make explicit the underlying bundle. 

The subspace of $\mathscr{A}$ of irreducible connections is, as usual, denoted by $\mathscr{A}^*$. There are slices for the action of $\G$ on $\mathscr{A}^*$ giving $\mathscr{B}^* = \mathscr{A}^*/\G$ the structure of a Banach manifold, see for instance \cite[Proposition 4.2.9]{DK} and the following discussion, or the slice theorem in \cite[Section 3]{FU}.

\subsection{The moduli space}\label{basic definition of the moduli space}
We denote by $F_A$ the curvature of the connection $A \in \mathscr{A}(\ad(P_E)) = \A(\su(E))$ and $F_A^+$ its self-dual part. Recall that it is equivariant with respect to the action of the gauge group $\G$. The anti-self-duality equation for $A \in \A$ is
\begin{equation}\label{asd eqn}
F_A^+ = 0 .
\end{equation}
The moduli space of anti-self-dual connections in $\ad(P)$ is defined to be the space 
\begin{equation*}
	M_{\kappa}^w := \left\{ [A] \in \mathscr{B}^w_\kappa | F_A^+ = 0 \right\}.
\end{equation*} 
Elements of this moduli space or particular representatives are referred to as {\em instantons}.
\\

There is an elliptic deformation complex associated to an instanton $[A]$:
\begin{equation*}
0 \to L^2_{l+1}(X;\su(E)) \xrightarrow{-d_A}
L^2_l(X;\Lambda^1 \tensor \su(E)) \xrightarrow{d^+_A}
L^2_{l-1}(X;\Lambda^2_+ \tensor \su(E)) \to
0 \ .
\end{equation*}
Here the second map is the derivative of the gauge-group action on $A$ and the third is the derivative of the map $\A \to L^2_{l-1}(X;\su(E))$ given by the left hand side of the equation (\ref{asd eqn}). The cohomology groups of this complex are, as usual, denoted by $H_A^0, H_A^1$ and $H_A^2$. An instanton $A$ is irreducible if and only if $H_A^0$ is zero. It is called {\em regular} if $H^2_A$ vanishes, and the moduli space $M_{\kappa}^{w}$ is regular if each instanton in it is regular. If an instanton $[A]$ is both irreducible and regular, then the moduli space $M$ has the structure of a smooth manifold in a neighbourhood of $[A]$. Its dimension is given by minus the index of the above complex or by the index of the `deformation operator' 
\begin{equation}\label{asd operator}
\delta_A := -d_A^* \oplus d_A^+ : L^2_l(X;\Lambda^1 \tensor \su(E)) \to L^2_{l-1}(X;\su(E)) \oplus \Lambda^2_+ \tensor \su(E)) \, .
\end{equation}
This index is given by the number
\begin{equation}\label{ex-dim instantons}
\begin{split}
d & := 8 \kappa + 3 \, (b_1(X) - b_2^+(X) -1) \\
	& = -2 \, \langle p_1(\su(E)) ,[X] \rangle +  3 \, (b_1(X) - b_2^+(X) -1) \, .
\end{split}
\end{equation}
This number is the {\em expected dimension} of the moduli space $M_\kappa^w$.
\\

In general the moduli space $M_{\kappa}^{w}$ need not be compact. However, there is the natural `Uhlenbeck-compactification' of it. In fact sequences of instantons can have sub-sequences whose curvatures become more and more concentrated at points of the manifold - the phenomenon called {\em bubbling}. As usual, we define an ideal instanton of instanton number $\kappa$ to be a pair ($[A],\bf{x})$, where $\bf{x}$ is an element of the n-fold symmetric product $\text{\em Sym}^n(X) = X^n/S_n$ (an unordered n-tuple of points in $X$) for some $n \geq 0$ and $[A]$ is an instanton in the moduli space $M_{\kappa-n}^w$. There is a topology on the space of ideal instantons 
\begin{equation*}
	{\text{I}}M_{\kappa}^{w} := \bigcup_{n=0}^{\infty} M^{w}_{\kappa-n} \times \text{ Sym}^n(X)
\end{equation*}
which is compact. In this topology each stratum admits its previously defined topology, but different strata are related by the notion of `weak convergence' \cite[\S 4.4. and in particular Condition 4.4.2]{DK}. The closure $\overline{M}^{w}_{\kappa} \subseteq \text{I}M_{\kappa}^{w}$ is therefore compact, see \cite[Theorem 4.4.3]{DK}.

\subsection{Holonomy perturbations}\label{Holonomy perturbations}
As the instanton moduli space we are studying consists of flat connections we cannot achieve transversality by perturbing the metric. A convenient choice of perturbation for this situation consists of holonomy perturbations, as used for instance in \cite{D,RS,Ta,K}. This is also the approach we shall choose, and so we will follow Kronheimer's exposition closely - here and in the section on $PU(2)$ monopoles.

Let $B$ be a closed 4-ball in $X$. Suppose a submersion
$q: S^1 \times B \to X$ which satisfies
 \begin{equation*}
   q(1,\_) = id_B \ 
 \end{equation*}
is given. 
Therefore each map $q_x:S^1 \to X$, defined by $q_x(z)=q(z,x)$, parametrises
a path in $X$ which is centered at $x$. Given a smooth connection $A \in \mathscr{A}$ the expression $\Hol_{q_x} (A)$ denotes the holonomy of the connection $A$ around the loop $q_x$. Therefore $\Hol_{q_x}(A)$ is an element of the group $\Aut(E)_x$. By letting $x \in B$ vary a unitary automorphism $\Hol_q(A)$ of $E$ over $B$ is given, i.e. a section of the  bundle $U(E)$ over $B$. The bundle $U(E)$ will be considered as a subbundle of the vector bundle $\gl(E)$ . Now let $\omega \in \Omega^2_+(X;\C)$ be a complex-valued self-dual two-form with compact support inside the ball $B$. Tensoring it with the section $\Hol(A)$ of $\gl(E)$ and extending it by zero onto $X$ a section
\[
\omega \tensor \Hol_q(A) \in \Omega^2_+(\gl(E)) \ , 
\]
is given. It defines, after applying the projection $\pi: \gl(E) \to \su(E)$, a section
\[
	V_{q,\omega}(A) \in \Omega^2_+(\su(E)) \ .
\]
There is an extension of this map to our configuration space of connections of class $L^2_l$ which admits uniform bound on its derivatives, once a reference connection $A_0 \in \A$ is fixed. Again, we recall that our considerations will require $l \geq 3$, although this is not a necessary condition for the following Proposition, see Remark \ref{conditions sobolev} below.

\begin{prop} \label{holonomy map}
For fixed $q$ and $\omega$ the map $V_{q,\omega}$ extends to a smooth map of
Banach manifolds
\begin{equation*}
 V_{q,\omega} : \mathscr{A}_l \to L^2_l(X,\Lambda^2_+\tensor \su(E))
 \ .
\end{equation*}
Furthermore there are uniform bounds on the derivatives of this map: There are constants $K_n$, depending only on $q$ and $A_0$, so that the n-th derivative
\[
  D^n V_{q,\omega} |_A: L^2_l(X,\Lambda^1 \tensor \su(E))^n \to
 L^2_l(X,\Lambda^2_+ \tensor \su(E))
\] 
satisfies
\[
 \norm{D^n V_{q,\omega} |_A (a_1, \dots, a_n)_{L^2_{l,A_0}}} \leq
   K_n \norm{\omega}_{C^l} \prod_{i=1}^{n} \norm{a_i}_{L^2_{l,A_0}} \ .
 \]
\end{prop}

Eventually a collection of submersions $q_i:S^1 \times B_i \to X$, $i \in \N $ 
as above will be chosen. Let $K_{n,i}$ be corresponding constants as guaranteed by the proposition above. Let $C_i$ be a sequence of numbers so that
 \[
   C_i \geq \sup\{ K_{n,i} | 1 \leq n \leq i \} \ .
 \]
Assume that $\omega_i$ is a sequence of self-dual two-forms with each $\omega_i$ having support inside $B_i$, and that the series
\[
  \sum_{i} C_i \norm{\omega_i}_{C^l}
\]
is convergent. Topologise the space of maps from $\mathscr{A}$ to
$L^2_l(X,\Lambda^2_+\tensor \su(E))$ with the $C^n$ semi-norms on bounded
subsets. It follows from the above proposition and the choice of the constants $C_i$ that the series 
\begin{equation}\label{convergence pert}
   \sum_i V_{q_i,\omega_i} 
\end{equation}
converges in the $C^n$-topology for each $n \in \N$. Therefore the series defines a smooth map of Banach manifolds
\[
 V_{{\boldsymbol{\omega}}}: \mathscr{A} \to L^2_l(X,\Lambda^2_+\tensor \su(E) \ . 
\]
We define the perturbation space correspondingly:
\begin{definition}
Fix the maps $q_i$, and the constants $C_i$ as above. The space $\mathscr{W}$ is defined to be the Banach space consisting of all sequences ${\boldsymbol{\omega}} = (\omega_i)_{i \in \N}$, with $\omega_i$ an element of the Banach space $\Omega^2_+(B_i)$,
so that the sum 
   \[
     \sum_i C_i \norm{\omega_i}_{C^l}
   \]
converges. We then define $V_{{\boldsymbol{\omega}}}$ to be the series $\sum V_{q_i,\omega_i}$.
\end{definition}


The dependence of $V_{\boldsymbol{\omega}}$ on $\bfomega$ is linear, and the map $\W \times \A \to L^2_l(X;\Lambda^2_+ \tensor \su(E))$, $(\bfomega,A) \mapsto V_{\bfomega}(A)$ is a smooth map of Banach-manifolds.
\\

For the following Remark see \cite[Section 3.2 and in particular the bottom of p. 70]{K}.
\begin{remark}\label{conditions sobolev}
Proposition \ref{holonomy map} continues to hold in the completions with the $L^{p}_{k}$ norms, for any $p \geq 1$ and $k \geq 0$. By a diagonalisation argument the constants $C_{i}$ in the preceding definition may be chosen so that the convergence (\ref{convergence pert}) occurs for all $p \geq 1$ and all $k \leq l$ at once, and we suppose in the sequel that such a choice has been made.
\end{remark}

Given $\bfomega \in \W$ the perturbed anti-self-duality equation for $A \in \A$ is now
\begin{equation}\label{perturbed asd}
 F_A^+ + V_{\boldsymbol{\omega}} (A) = 0 \ .
\end{equation}
Correspondingly we define the moduli space perturbed by $\bfomega$ to be the space
\begin{equation*}
 M^w_{\kappa}({\boldsymbol{\omega}}) := \{ [A] \in B^w_\kappa | F_A^+ + V_{\boldsymbol{\omega}}(A) = 0 \} \ .
\end{equation*}
\\

There is an elliptic deformation complex associated to an instanton $[A]$ in the perturbed moduli space $M_{\kappa}^{w}(\bfomega)$ also:
\begin{equation*}
0 \to L^2_{l+1}(X;\su(E)) \xrightarrow{-d_A}
L^2_l(X;\Lambda^1 \tensor \su(E)) \xrightarrow{d^+_{A,\bfomega}}
L^2_{l-1}(X;\Lambda^2_+ \tensor \su(E)) \to
0 \ ,
\end{equation*}
where now $d_{A,\bfomega}^+ = d_A^+ + d V_{\bfomega}|_A$ is the derivative of the map $\A \to L^2_{l-1}(X;\su(E))$, $A \mapsto F^+_{A} + V_{\bfomega}(A)$ at the instanton $A$. The deformation operator of the $\bfomega$-perturbed equations is $\delta_{A,\bfomega}:= -d_A^* \oplus d^+_{A,\bfomega}$. 

Note that $d^+_{A,\bfomega}$ differs from 
$d^+_{A}$ only by the addition of a compact operator $L^2_l(X;\Lambda^1 \tensor \su(E)) \rightarrow L^2_{l-1}(X;\Lambda^2_+ \tensor \su(E))$; the index of this elliptic complex is therefore the same as that of the unperturbed anti-self-duality equations above. 
The cohomology groups of this complex are denoted by $H_A^0, H_{A,\bfomega}^1$ and $H_{A,\bfomega}^2$. Again, an instanton is called {\em regular} if the cohomology group $H^2_{A,\bfomega}$ vanishes. Likewise, local models for the moduli space $M(\bfomega)$ show that in a neighbourhood of a point $[A]$ which is both irreducible and regular the moduli space admits the structure of a smooth manifold of the expected dimension.

\subsection{Compactness}
The existence of an Uhlenbeck-type compactification of the perturbed moduli space $M^{w}_{\kappa}(\bfomega)$ was proved in \cite[Proposition 3.5]{K}:

\begin{proposition}\label{compactness instantons}
	Let $A_n$ be a sequence of connections in the Hermitian bundle $E \to X$ representing points $[A_n]$ in the moduli space $M^{w}_{\kappa}(\bfomega)$. 
	Then there is a point ${\bf x} \in {\text{\em Sym}}^n(X)$, a connection $A'$ in a bundle $E' \to X$ representing an element of a moduli space $M^{w}_{\kappa-n}(\bfomega)$ with the following property: After possibly passing to a subsequence,
\begin{enumerate}
 
\item[(i)] there are bundle isomorphisms
\begin{equation*}
 	h_n: E|_{X \setminus {\bf x}} \to E'|_{X \setminus {\bf x}} \ ,
\end{equation*}
so that $(h_n)_*(A_n)$ converges to $A'$ in $L^p_1(K)$ for all compact subset $K \subseteq X \setminus {\bf x}$, and for all $p \geq 2$, and
\\
\item[(ii)] the sequence of measures $(\abs{F_{A_n}}^2 vol_g)$ converges to the measure $\abs{F_{A'}}^2 \text{vol}_g + 8 \pi^2 \sum_{x \in {\bf x}} \delta(x)$ in the weak-*-topology of measures, i.e. for any continous function $f$ on $X$ one has 
\[
	\int_X f\,  \abs{F_{A_n}}^2 \text{vol}_g \ \to \int_X f\, \abs{F_{A'}}^2 \text{vol}_g + 8 \pi^2 \sum_{x \in {\bf x}} f(x) \ 
\]
as $(n \to \infty)$. 
\end{enumerate}
\end{proposition}
That we have a weaker notion of convergence than in \cite[\S 4.4]{DK} is due to the fact that changes of the connection $A$ have effect globally on the section $V_{\bfomega}(A)$. 
\\

An ideal instanton of instanton number $\kappa$ is a pair $([A],\bf{x})$ where ${\bf x} \in {\text{\em Sym}}^n(X)$ and $[A]$ is an element of the perturbed moduli space $M^w_{\kappa-n}(\boldsymbol{\omega})$. The space of $\bfomega$-perturbed ideal instantons of instanton number $\kappa$ is defined to be
\begin{equation*}
	{\text{I}}M_{\kappa}^{w}(\bfomega) := \bigcup_{n=0}^{\infty} M^{w}_{\kappa-n}(\bfomega) \times \text{Sym}^n(X) \ ,
\end{equation*}
with the notion of convergence between the different strata as in Proposition \ref{compactness instantons}, and with each stratum having its original topology. It follows from the above Proposition that the closure of the moduli space $M^{w}_{\kappa}(\bfomega)$ inside the space of ideal instantons ${\text{I}}M_{\kappa}^{w}(\bfomega)$ is compact. 

\subsection{Transversality}
Under an additional condition on the set of submersions $q_i: S^1 \times B_i \to X$ transversality can be achieved for the moduli space. The condition is as follows:
\begin{condition} \label{condition}
For any point $x \in X$ the set of loops
\begin{equation*}
 	\left\{ q_i|_{S^1 \times \{x\}} | i \in \N, \, x \in \text{int}(B_i) \right\}
\end{equation*}
is $C^1$-dense in the space of smooth loops based at $x$. 
\end{condition}
It is a non-trivial exercise to convince oneself that this condition can always be fulfilled. This is implicit in \cite{K}. 
\begin{theorem}\label{transversality instantons} \cite{K}
Suppose the submersions $q_i: S^1 \times B_i \to X$ satisfy the above condition. Then the smooth map of Banach manifolds
\begin{equation*}
\begin{split}
	g: \W \times \A^* & \to L^2_{l-1}(X;\Lambda^2_+ \tensor \su(E)) \\
	({\boldsymbol{\omega}}, A) & \mapsto F_A^+ + V_{\boldsymbol{\omega}}(A)
\end{split}
\end{equation*}
is transverse to zero. 
\end{theorem}
The key-point in the proof is that for an irreducible connection $A \in \A^*$ and point $x \in X$ the holonomy-sections $\Hol_{q_i}(A)$, associated to submersions $q_i:S^1 \times B_i \to X$ so that $x$ is contained in the interiour of $B_i$, span $\gl(E)_x$. Furthermore, after exhibiting a basis out of these sections, this basis continues to be a basis of $\gl(E)$ in a neighbourhood of $x$. It is at this point that the inclusion $L^2_l \hookrightarrow C^0$ is used. \\

Let us denote by $\mathscr{M} := g^{-1}(0) /\G$ the {\em parametrised moduli space}. Applying the Sard-Smale theorem to the projection $\mathscr{M} \to \W$ yields the following result in the standard way:

\begin{corollary}\label{generic regularity}
For a residual set of perturbations $\bfomega \in \W$ the moduli space 
$M_{\kappa}^{w,*}(\bfomega)$ is regular for all $w\in H^2(X;\Z)$ and instanton numbers $\kappa$. It therefore admits the structure of a smooth manifold of the expected dimension d given by (\ref{ex-dim instantons}).  
\end{corollary}

As usual, a `residual' subset of a complete metric space is a countable intersection of open and dense sets. By Baire's theorem such a set is dense itself. 
 
\subsection{Orientations}\label{orientations instantons}
As in Donaldson's first applications of gauge theory to 4-manifold topology \cite{D2} the moduli space is given an orientation by a choice of orientation for a real determinant line bundle of a family of Fredholm operators. 
\\

The determinant line of a Fredholm operator $T: V \to W$ is given by $\det(\ker(T)) \tensor \det(\coker(T))^*$, where $\det(F)$ denots the maximal exteriour power of a finite dimensional vector space $F$, and $\det(F)^*$ its dual. For a family of Fredholm operators $(T_c:V \to W)_{c \in C}$, parametrised continously by a topological space $C$, there is a line bundle $\det(T)$ on the space $C$ whose fibre over the point $c \in C$ is given by the line $\det(T_c)$, and whose topology is given as described in the next paragraph.

If $T_{c_0}$ has trivial cokernel then the family of kernels $\ker(T_c)$ admits the structure of a vector bundle in a natural way over a neighbourhood of $c_0$, so $\det(T)$ admits a natural topology when restricted to that neighbourhood. If $T_{c_0}$ has non-trivial cokernel we pick a subspace $J \subseteq W$ that surjects onto the cokernel of $T_{c_0}$. The space $J$ also surjects onto the cokernels of $T_c$ for $c$ out of a neighbourhood of $c_0$, and there is a natural exact sequence 
\begin{equation*}
 0 \to \ker T_c \to T_c^{-1}(J) \to J \to \coker T_c \to 0 \ ,
\end{equation*}
where the third of the five maps involved is induced by $T_c$. It is an algebraic fact (see for instance \cite[Section 20.2]{KM_floer} or \cite{DK}) that this implies the existence of a natural isomorphism 
\begin{equation*}
 \det(\ker(T_c)) \tensor \det(\coker(T_c))^{*} \cong \det(T_c^{-1}(J)) \tensor \det(J)^{*} \ .
\end{equation*}
Now in a neighbourhood of $c_0$ the family $T_c^{-1}(J)$ forms naturally a vector bundle, and so there is a natural structure of real line bundle on $\det(T_c^{-1}(J)) \tensor \det(J)^{*}$, over this neighbourhood. It is possible to {\em define} the topology of $\det(T)$ on this neighbourhood by that of $\det(T_c^{-1}(J)) \tensor \det(J)^{*}$. Indeed, any choice $K \subseteq W$ with $J \subseteq K$ yields the same topology by the above construction. One can also see that there are then continous transition functions on overlaps. We refer the reader to \cite[Section 20.2]{KM_floer} and \cite[Appendix A]{Fro3} for a detailed discussion.
\\

In our situation, the determinant line bundle formed by the family of Fredholm operators 
\begin{equation*}
	\delta_{A,\bfomega}= -d_A^* \oplus d_{A,\bfomega}^+ : L^2_{l}(X;\Lambda^1 \tensor \su(E)) \to L^2_{l-1}(X;\su(E) \oplus \Lambda^2_+ \tensor \su(E)) \ 
\end{equation*}
 defined
 on the space of connections $\A$ is relevant for orientations. 
The restriction of the determinant line bundle of this family of operators to $\A^*$ 
 descends to the quotient $\bonf^*$. We denote by $\Lambda_{\bfomega}$ this line bundle. Its restriction to the regular part of the moduli space $M^*(\bfomega)$ is equal to its orientation line bundle. In fact, at these points, the cokernels of the deformation operator vanish and the kernels $H^1_{A,\bfomega}$ are precisely the tangent spaces to the moduli space $M^{*}(\bfomega)$. In particular, a regular moduli space $M^{*}(\bfomega)$ is orientable if $\Lambda_{\bfomega}$ is orientable, and a trivialisation of $\Lambda_{\bfomega}$ provides a preferred orientation for the moduli space $M^*(\bfomega)$. 

Now the space of perturbations $\W$ is contractible, so the line bundles $\Lambda_{\bfomega}$ on $\bonf^*$ corresponding to different perturbations $\bfomega \in \W$ are canonically isomorphic. Donaldson's Theorem \cite[Corollar 3.27]{D} or \cite[Proposition 5.4.3]{DK} states that the line bundle $\Lambda:=\Lambda_0 \to \bonf^{* w}_{\kappa}$ is indeed trivial, and an orientation is determined by a `homology orientation' of $X$ (that is, an orientation of the real vector space $\mathscr{H}^1(X;\R) \oplus \mathscr{H}^2_+(X;\R)$ of harmonic one-forms and self-dual two-forms), see also  \cite[Section 5.4]{DK}. We denote by the letter $o$ a choice of trivialisation of the line bundle $\Lambda$. 

\subsection{Moduli spaces over negative definite four-manifolds} \label{moduli space on neg def manifolds}

We restrict now our attention to smooth, closed Riemannian four-manifolds $X$ with $b_2^+(X) = 0$ and $b_2(X) \geq 1$. According to the theorem of Donaldson's \cite{D} the intersection form of such a four-manifold is equivalent to the diagonal one. Kim Fr\o yshov pointed out to the author that there is a generalisation of the condition in \cite[Section 4.2.1]{T} to assure the absence of reducibles in the moduli space: 
\begin{lemma}\label{no reductions}
 Suppose $E \to X$ is a Hermitian bundle of rank 2 so that the square of its first Chern class $\langle c_1(E)^2,[X] \rangle$ is divisible by 4, so that the second Chern class satisfies 
\begin{equation*}
 \langle c_2(E), [X] \rangle = \frac{1}{4} \, \langle c_1(E)^2,[X] \rangle \ - k \, ,
\end{equation*} 
for some integer $k \geq 0$, 
and so that in $H^2(X;\Z)/\text{\em Torsion}$ the class $c_1(E)$ is not divisible by 2. Then $E$ does not admit a topological decomposition $E = L \oplus K$ into the sum of two complex line bundles $L$ and $K$. 
\end{lemma}
{\em Proof:}
Suppose we have a decomposition into line bundles $E = L \oplus K$. Let $z_1 := c_1(L)$ and $z_2:= c_1(K)$ be the first Chern-classes of the line bundles. Notice that by the assumption we have
\begin{equation}
\begin{split}
 \langle (z_1 - z_2)^2,[X] \rangle & = \langle (z_1 + z_2)^2 - 4 \; z_1 \, z_2 ,[X] \rangle \\
	& = \langle c_1(E)^2 - 4 \, c_2(E), [X] \rangle \\
	& = 4 \,  k \ .
\end{split}
\end{equation}
If $k=0$ this equation and the fact that the intersection form is definite show that  $z_1 - z_2$ is a torsion class. But $c_1(E) = z_1 + z_2 $, contradicting the non-divisibility condition. If $k > 0$ this gives right off a contradiction to the fact that $X$ is negative definite. 
\qed
%
\begin{corollary}
 Let $E \to X$ be as in the previous lemma. Then the associated moduli space $M^w_0$ does not admit reducibles. The same holds for the perturbed moduli spaces $M^w_{0}(\bfomega)$.  
\end{corollary}
\begin{observation}
A bundle $E \to X$ with the characteristic classes as in the preceding lemma exists if and only if $b_2(X) \geq 4$.  
\end{observation}

For a connection $A \in \mathscr{A}(\su(E))$ Chern-Weil theory gives the following formula:
\begin{equation}\label{chern-weil formula}
\begin{split}
\frac{1}{8\pi^2} (\norm{F_A^-}^2_{L^2(X)} - \norm{F_A^+}^2_{L^2(X)})   & = \, - \frac{1}{4} \,\langle  p_1(\su(E)) , [X] \rangle \\
	& = \langle c_2(E) - \frac{1}{4} c_1(E)^2 , [X] \rangle
\end{split}
\end{equation}
In particular, for anti-self-dual connections the left hand side of this equation is always non-negative. 

\begin{proposition}\label{our moduli space}
Suppose the negative definite four-manifold $X$ has first Betti-number $b_1(X)=1$ and admits a class $w \in H^2(X;\Z)$ so that $\langle w^2,[X]\rangle$ is divisible by four, and so that (modulo torsion) $w$ is not divisible by 2. Then the moduli space $M^w_0$, associated to the bundle $E \to X$ with $c_1(E) = w$ and $\langle c_2(E),[X]\rangle = 1/4 \,  \langle w^2,[X] \rangle$, is compact, does not contain reducibles, consists of flat connections in $\su(E)$, and is of expected dimension zero. 
\end{proposition}
{\em Proof:} Equation (\ref{chern-weil formula}) implies that all the lower strata of the Uhlenbeck-com\-pac\-ti\-fi\-cation of $M_{0}^{w}$ are empty, so $M_{0}^{w}$ must already be compact. The remaining claims follow from the above Lemma \ref{no reductions} and the dimension-formula (\ref{ex-dim instantons}). \qed
\begin{remark}
 The `Casson-type' instanton moduli spaces appearing in \cite{T} are as in Proposition \ref{our moduli space}, but associated to elements $w = \sum e_i$, where the $\{ e_i \}$ yield a basis of $H^2(X;\Z)/\text{\em Torsion}$ diagonalising the intersection form. This requires then $b_2(X) \equiv 0 \ (\text{mod} \ 4)$. 
\end{remark}

\subsection{Definition of the invariant} \label{definition casson}
We shall stay in the situation of the preceding section, and in particular in that of Proposition \ref{our moduli space}. If the compact moduli space $M^{w}_{0}$ were regular (and zero-dimensional) we would define an integer by a signed count of its finite number of elements (each regular point is isolated). In general, we will have to consider perturbations of this moduli space: 

\begin{proposition}\label{our moduli space perturbed}
Suppose the $C^0$-norm of the perturbation $\bfomega \in \W$ is sufficiently small. Then the moduli space $M_{0}^{w}(\bfomega)$ is compact. If further $\bfomega$ is chosen among a residual subset of $\W$ so that the conclusion of Corollary \ref{generic regularity} holds, then $M_{0}^{w}(\bfomega)$ is a compact zero-dimensional manifold.
\end{proposition}
{\em Proof:} The claim on compactness is an easy consequence of the Chern-Weil formula (\ref{chern-weil formula}) and the structure of the compactification of the moduli space, Proposition \ref{compactness instantons}.
\qed

The following is the fundamental definition of this paper, so we rephrase all conditions we have imposed so far and state some of the consequences already established.

\begin{definition} \label{casson invariant perturbation}
Suppose the negative definite four-manifold $X$ has first Betti-number $b_1(X)=1$ and admits a class $w \in H^2(X;\Z)$ so that $\langle w^2,[X]\rangle$ is divisible by four, and so that (modulo torsion) $w$ is not divisible by 2 (this condition requires $b_2(X) \geq 4$.) The moduli space $M^w_0$, associated to the bundle $E \to X$ with $c_1(E) = w$ and $\langle c_2(E),[X]\rangle = 1/4 \,  \langle w^2,[X] \rangle$, is compact, does not contain reducibles, consists of flat connections in $\su(E)$, and is of expected dimension zero (by Proposition \ref{our moduli space}).

Let $\bfomega \in \W$ be a perturbation so that the conclusion of Proposition \ref{our moduli space perturbed} holds. In particular, the moduli space $M_0^w({\bfomega})$ is a compact, zero-dimensional manifold.

Suppose an orientation $o$ of the determinant line bundle $\Lambda \to \bonf^*$ is chosen, and, therefore an orientation for the moduli space $M_{0}^{w}(\bfomega)$, according to Section \ref{orientations instantons}. Then we define the number $n_{o}(\bfomega)$ as the signed count of the moduli space $M_{0}^{w}(\bfomega)$,
\begin{equation*}
	n_{w,o}(\bfomega) := \# M_{0}^{w}(\bfomega) \ .
\end{equation*}
Here an instanton $[A] \in M_{0}^{w}(\bfomega)$ is counted with $+1$ if the orientation of the determinant line $\det(\delta_{A,\bfomega})$ at $[A]$ determined by $o$ coincides with the preferred orientation 
\begin{equation*}
\det(\delta_{A,\bfomega}) = \det(\ker(\delta_{A,\bfomega})) \tensor \det(\coker(\delta_{A,\bfomega}))^* = \R
\end{equation*}
determined by the trivial kernel and cokernel of $\delta_{A,\bfomega}$ at the irreducible and regular point $[A]$, and with $-1$ in the opposite case.
\end{definition}

It is worth noting (see for instance the Appendix A of \cite{S}) that the relative sign between two instantons $[A_0]$ and $[A_1]$ can be computed from the number of crossings $\mu=\sum_t \dim(\ker(\delta_{A_t,\bfomega}))$ of a generic path $t \mapsto A_t$. It is $(-1)^{\mu}$. 
\\

The definition of the perturbation space $\W$ depends on the choice of a Riemannian metric $g$ on $X$. 
In the following proposition we shall write $\W_g$ to signify this dependence. Furthermore, the Riemannian metric also goes into the definition of the moduli space $M^w_0(\bfomega)$ via the perturbed anti-self-duality equations, depending on the Riemannian metric. To express this we shall write $n_o(\bfomega,g)$ for the number defined above. 

\begin{proposition}\label{well-definedness} 
Let $g, g'$  be Riemannian metrics on $X$ and suppose perturbations $\bfomega \in \W_g$, $\bfomega' \in \W_{g'}$ are chosen so that the conclusion of Proposition \ref{our moduli space perturbed} holds in each case. Then we have
\begin{equation*}
	n_{w,o} (\bfomega,g) = n_{w,o} (\bfomega',g') \ .
\end{equation*}
\end{proposition}
{\em Proof:}
The proof follows from the standard cobordism argument in such kind of situations. See for instance \cite[page 23]{K} or \cite[section 7.1]{T5}. \qed

\begin{definition} \label{casson invariant without perturbation}
The number $n_{w,o}(\bfomega,g)$ is therefore independent of the underlying Riemannian metric on $X$ and the chosen perturbation $\bfomega$. It only depends on the topology of the smooth manifold $X$ and the cohomology class $w$ (determining the bundle $E \to X$ which defines the moduli space) and a choice of orientation $o$. It is therefore convenient to denote this number by $n_{w,o}(X)$. 
\end{definition}

\begin{proposition}\label{tensor product with line bundle}
Let $c \in H^2(X;\Z)$ be an element in the second cohomology group. 
Suppose $o$ is a choice of trivialisation of the bundle $\Lambda_0 \to \bonf^{w,*}_{0}$, and suppose $o'$ is
a choice of trivialisation of the bundle $\Lambda_0' \to \bonf^{w+2c,*}_{0}$. Then we have
\begin{equation*}
 n_{w,o}(X) \, = \, \pm \, n_{w+2c,o'}(X) \ .
\end{equation*}
\end{proposition}
{\em Proof:}
Let $L \to X$ be a Hermitian line bundle with $c$ as its first Chern class. Let $a_0$ be a Hermitian
connection in $L$. The map
\begin{equation*}
 \begin{split}
  	\A_{a}(E) & \to \A_{a \tensor a_0^2}(E \tensor L) \\
	A & \mapsto A \tensor a_0
 \end{split}
\end{equation*}
is an isomorphism of affine spaces. It descends to a homeomorphism $\bonf^{w}_{\kappa} \to
\bonf^{w+2c}_{\kappa}$ and induces a homeomorphism of moduli spaces
$M^{w,*}_{\kappa}(\bfomega) \cong M^{w+2c,*}_{\kappa}(\bfomega)$. This is either
orientation-preserving or orientation-reversing. In our case where the index of the deformation operator is
zero this is easily seen by the interpretation of the relative orientations in terms of crossing numbers. 
\qed

\subsection{Moduli spaces with cut off perturbations} \label{moduli spaces with cut off perturbations}
For reasons that will become apparent later we define the following modification of the moduli spaces of perturbed instantons.

Let $h : \R_+ \to [0,1]$ be a smooth real-valued function defined on the half-line of non-negative real numbers with the following properties:
\begin{enumerate}
	\item $ h(t) = 1$ for $0 \leq t \leq 1$,
	\item $ h(t) = 1/t$ for $t \geq 2$. 
\end{enumerate}

Let $m > 2$. The equation we shall now consider for 
 $\bfomega \in \W$ and $A \in \A$ is 
the following: 
\begin{equation}\label{perturbed asd cut off}
 F_A^+ + h(\norm{F_{A}}_{L^m(X)}) \, V_{\boldsymbol{\omega}} (A) = 0 \ .
\end{equation}

Correspondingly, we define the moduli space perturbed by $\bfomega$ to be the space
\begin{equation} \label{perturbed asd cut off}
 M^w_{\kappa}({\boldsymbol{\omega}},h,m) := \{ [A] \in B^w_\kappa | F_A^+ + h(\norm{F_{A}}_{L^m(X)}) \, V_{\boldsymbol{\omega}}(A) = 0 \} \ .
\end{equation}

Our previous discussion on compactification of this moduli space has one significant difference. Namely, in the lower strata the perturbation term vanishes as a consequence of the following Lemma:

\begin{lemma}\label{L^m bound on curvature}
Suppose the sequence of measures $\abs{F_{A_n}}^2 \text{vol}_g$ converges weakly to a measure $\nu$ on $X$, in the weak-*-topology of measures: For any continous function $f$ on $X$ one has
\begin{equation*}
	\int_X f\,  \abs{F_{A_n}}^2 \text{vol}_g \ \to \int_X f\,  \nu \ \ 
\end{equation*}
as $n \to \infty$. 
Suppose the point $x \in X$ is so that it is not contained in a geodesic ball  of $\nu$-measure less than $\epsilon^2$, for some $\epsilon > 0$. Then for any $m > 2$ the sequence of norms $\norm{F_{A_n}}_{L^m(X)}$ tends to infinity.
\end{lemma}
{\em Proof:}
Suppose this were not the case. Then there would be a constant $N > 0$ and  a subsequence $(n_k)$ of $(n)$ so that 
\begin{equation*}
\norm{F_{A_{n_k}}}_{L^m(X)} \leq N 
\end{equation*}
held for all $k \in \N$. Let $B(x,r)$ be the geodesic ball around $x$ of radius $r > 0$. By the H\"older inequality we get 
\begin{equation*}
\begin{split}
	\norm{F_{A_{n_k}}}^2_{L^2(B(x,r))} 
	\leq \text{vol}_g (B(x,r))^{\frac{m-2}{m}} \ \norm{F_{A_{n_k}}}^2_{L^m(B(x,r))} \ .
\end{split}
\end{equation*}
Clearly, the expression on the left-hand side converges to  the $\nu$ measure of the ball $B(x,r)$. On the other hand, by our assumption on the $L^m$ boundedness, the expression on the right-hand side tends to zero as the radius $r$ tends to zero. But this contradicts the hypothesis on the point $x$. \qed

There are therefore obvious modifications of Proposition \ref{compactness instantons} and Corollary \ref{generic regularity}. However, compactness of the Casson-type moduli space $M^{w}_{0}(\bfomega)$ follows from the Chern-Weil formula as long as the $L^\infty$ norm of the perturbation $\bfomega$ is small enough, and therefore the results of section \ref{moduli space on neg def manifolds} continue to hold for the moduli space $M^{w}_{0}(\bfomega,h,m)$. Correspondingly, there is an invariant $n_{w,o}(\bfomega,h,m)$, and the following Proposition follows from a standard cobordism argument again.
\begin{proposition}
For the two numbers $n_{w,o}(\bfomega,h,m)$ and $n_{w,o}(\bfomega)$ the following equation holds:
\[
n_{w,o}(\bfomega,h,m) = n_{w,o}(\bfomega) \ .
\]
\end{proposition}

\section{Moduli spaces of $PU(2)$ Seiberg-Witten monopoles} 
Here we shall recall the $PU(2)$-monopole equations and their moduli space associated to the data of a $Spin^c$-structure $\mathfrak{s}$ and a Hermitian bundle $E \to X$ of rank $2$ on a Riemannian four-manifold $X$. We shall define the configuration space and the moduli space and recall how to get a uniform bound on the spinor component of a solution to the monopole equations. We then show how the equations are perturbed, sketch the Uhlenbeck compactification for the perturbed moduli space and show how to obtain transversality. Furthermore, we shall show how to define a preferred orientation on the irreducible part of the moduli space. At least in slightly different situations these results are already well-known \cite{FL}, \cite{T2}. 


\subsection{The configuration space}
Let $X$ be a closed oriented Riemannian four-manifold with a
$Spin^c$ structure $\mathfrak{s}$ on it. The $Spin^c$ structure consists of two Hermitian rank 2 vector bundles $S^\pm_{\mathfrak{s}}$ with identified determinant line bundles and a Clifford multiplication 
\begin{equation*} 
\gamma : \Lambda^1(T^*X) \to \Hom_\C(S^+_\mathfrak{s},S^-_\mathfrak{s}) \ , 
\end{equation*}
see for instance Witten's original work \cite{Witten} or Kronheimer and Mrowka's book \cite{KM_floer}. 

Let us furthermore suppose we are given a
Hermitian vector bundle $E$ with determinant line bundle $w = \det(E)$ on $X$. We can then form the `twisted' spinor bundles 
\begin{equation*}
W^\pm_{\data}:= S^\pm_{\mathfrak{s}} \tensor E .
\end{equation*}
Clifford multiplication extends by tensoring with the identity on $E$. 

We continue to denote by $\A$ the space of connections in $E$ inducing a fixed connection in the determinant line $w$ of class Sobolev class $L^2_{l}$. 
We define our pre-configuration space to be the product
\begin{equation*}
\mathscr{C} := \A \times L^2_{l}(X;S^+ \tensor E) \ .
\end{equation*}
As in the last section, we will require $l \geq 3$ ensuring an embedding $L^2_l \hookrightarrow C^0$. 

According to our above notation convention we use the notation $\mathscr{C}_{\kappa,\mathfrak{s}}^{w}$ if we want to emphasise that this configuration space is associated to the topological data consisting of a $Spin^c$ structure $\s$ and a bundle $E \to X$ with instanton number $\kappa$ as in the formula (\ref{instanton number}) and determinant line bundle $w:= \det(E)$. 

The gauge group $\G$ acts naturally on $\conf$: It acts on $\A$ as described in Section \ref{elementary definition configuration space} above, and by $(u,\Psi) \mapsto (\id \tensor u) (\Psi)$ on the spinor section. The quotient 
\[
	\mathscr{B}_{\kappa,\s}^{w} := \mathscr{C}_{\kappa,\s}^{w} / \G \ .
\]
is called the configuration space of $PU(2)$ monopoles. By $\Gamma_{(A,\Psi)} \subseteq \G$ we shall denote the stabiliser of the pair $(A,\Psi)$ under the gauge group action. We shall define by $\mathscr{C}^*$ and by $\mathscr{B}_{\s}^*$ the spaces corresponding to {\em trivial} stabiliser. Furthermore, we shall define by $\mathscr{C}^{**}$ the subspace of pairs $(A,\Psi)$ with $A$ irreducible, $\Gamma_A = \Z/2$, and with $\Psi \neq 0$. Certainly $\mathscr{C}^{**} \subseteq \mathscr{C}^*$ but the converse is not true.  

\subsection{A quadratic map on the spinor component}
The bundle $S^+ \tensor E$ is modelled on $\C^2_+ \otimes \C^2$.
We define the orthogonal projections
\begin{equation*}
\begin{split}
  P: \ \gl(\C^2 \tensor \C^2) & \to \sl(\C^2) \tensor \sl(\C^2) \
\end{split}
\end{equation*}
as being the tensor product of the two orthogonal projections $\gl(\C^2) \to \sl(\C^2)$. We then define the map
\begin{equation*}
\begin{split}
	\mu: \C^2_+ \tensor \C_2 \oplus \C^2_+ \tensor \C_2 &  \to \sl(\C^2) \tensor \sl(\C^2) \\
	(\Psi,\Phi) & \mapsto P( \Psi \Phi^*) \ ,
\end{split} 
\end{equation*}
where $(\Psi \Phi^*) \in \gl(\C^2_+ \tensor \C^2)$ is defined to be the
endomorphism $\Xi \mapsto
\Psi (\Phi,\Xi)$, where $(\Phi,\Xi)$ stands for the (standard) inner product of the two elements. Instead of $\mu(\Psi,\Psi)$ we shall just write $\mu(\Psi)$, where then $\mu$ is a quadratic map. Both the bilinear and the quadratic map are equivariant with respect to the structure groups of $S^+$ and $E$ and so define well-defined maps
\[
\mu : S^+ \tensor E \oplus S^+ \tensor E \to \sl(S^+) \tensor \sl(E) \ ,
\]
respectively 
\[
\mu : S^+ \tensor E \to \su(S^+) \tensor_\R \su(E) \ .
\]
Note that on the level of bundles this map is also equivariant with respect to the action of the gauge group $\G$. 

\subsection{The Dirac operator, moduli spaces of $PU(2)$ monopoles} \label{properties of monopoles}

We suppose $B$ is a fixed $Spin^c(4)$ connection on $S^+ \oplus S^-$ with respect to some $SO(4)$ connection on the tangent bundle $TX$, not necessarily equal to the Levi-Civita connection. As the particular choice of $Spin^c$ connection will be kept fixed in what follows it is suppressed from notation. The reason for this somewhat unusual choice will become apparent later when we are concerned with regularity in lower strata of the Uhlenbeck compactification. 

Composing the Clifford multiplication with the tensor product connection $\nabla_B \tensor \nabla_A$ yields the `twisted' Dirac operator
\begin{equation*}
 \dirac_{A}^+:= \gamma \circ (\nabla_B \tensor \nabla_A) : L^2_{l}(X;S^+ \tensor E) \to L^2_{l-1}(X;S^- \tensor E) \ . 
\end{equation*}

Let furthermore $\beta \in \Omega^1(X;\C)$ be a complex-valued 1-form, and let $K$ be a general linear orientation-preserving automorphism of the real rank-3 bundle $\su(S^+_\s)$. 

The $PU(2)$ monopole equations we consider, associated to a pair $(A,\Psi) \in \conf$, are given by
\begin{equation} \label{pu2 monopoles}
	\begin{split}
		\dirac^+_A \Psi + \gamma(\beta) \Psi & = 0 \\ \gamma ( F_A^+) - K (\mu(\Psi)) & = 0 \ .
	\end{split}
\end{equation}
Solutions of this equations will also be called `monopoles'.

The left-hand side of these equations can be considered\footnote{Here and later we will feel free to use the isomorphism $\gamma:\Lambda^2_+ \to \su(S^+)$ without making it explicit in the notation} as a map $g:\conf \to L^2_{l-1}(X;S^- \tensor E \oplus \Lambda^2_+(X) \tensor \su(E))$. This map is equivariant with respect to the action of the gauge group. The moduli space of $PU(2)$ monopoles is defined to be the space
\begin{equation} \label{monopoles}
	M_{\kappa,\s}^{w}:=\left\{ [A,\Psi] \in \bonf^{w}_{\kappa,\s} | (\ref{pu2 monopoles}) \text{ holds} \right\} / \G.
\end{equation}
Again, there is an elliptic deformation complex associated to a monopole $(A,\Psi)$:
\begin{equation*}
\begin{split}
0 \to L^2_{l+1}(X;\su(E)) & \xrightarrow{d^0_{(A,\Psi)}}  
L^2_l(X;S^+ \tensor E \oplus \Lambda^1 \tensor \, \su(E)) \\
& \xrightarrow{d^1_{(A,\Psi)}}
L^2_{l-1}(X;S^-\tensor E \oplus \Lambda^2_+ \tensor \, \su(E)) \to
0 \ .
\end{split}
\end{equation*}
Here $d^0_{(A,\Psi)}$ is the derivative of the gauge-group action $u \mapsto u(A,\Psi)$ at the identitiy, and $d^1_{(A,\Psi)} = dg|_{(A,\Psi)}$, the derivative of the monopole map $g$ at a solution $(A,\Psi)$. As in instanton theory, the kernel of the operator $d^{0,*}_{(A,\Psi)}$ provides a slice of the action of the gauge group $\G$ on $\conf$ in a neighbourhood of the orbit through $(A,\Psi)$. Again, a configuration $(A,\Psi)$ is called irreducible if the zeroth cohomology space of the above complex vanishes, $H^0_{(A,\Psi)}=0$, and regular if $H^2_{(A,\Psi)}=0$.

The deformation operator
\begin{equation*}
	D_{(A,\Psi)} = d^{0,*}_{(A,\Psi)} \oplus d^1_{(A,\Psi)}
\end{equation*}
is elliptic. Its kernel is given by the cohomology space $H^1_{(A,\Psi)}$ of the above complex and its cokernel is given by the sum of cohomology spaces $H^0_{(A,\Psi)} \oplus H^2_{(A,\Psi)}$. If this cokernel vanishes then the slice theorem together with implicit function theorems show that the moduli space has the structure of a smooth manifold in a neighbourhood of $[A,\Psi]$, of dimension given by the index of the elliptic operator. 

The deformation operator is homotopic to the operator $D_{(A,0)}$ which has the simple form
\begin{equation*}
	D_{(A,0)} = -d_A^* \oplus d_A^+ \oplus \dirac^+_A \, ,
\end{equation*}
as a map $L^2_{l}(X;\Lambda^1 \tensor \, \su(E) \oplus S^+ \tensor E) \to L^2_{l-1}(X;(\Lambda^0 \oplus \Lambda^2_+)\tensor \, \su(E) \oplus S^- \tensor E)$. It is equal to the sum of the instanton deformation operator $\delta_A = -d_A^* \oplus d_A$ and the Dirac operator $\dirac^+_A$. The expected dimension of the moduli space $M_{\kappa, \s}^{w}$ is therefore given by the index of this direct sum of two elliptic operators,
\begin{equation*}
	\text{ex-dim}(M_{\kappa,\s}^{w}) = \ind(\delta_A) \oplus 2 \, \ind_\C(\dirac^+_A) \ , 
\end{equation*}
where the operator $\delta_A$ has been computed in (\ref{ex-dim instantons}) above and the index of the Dirac operator will be computed in the next section.\\

We terminate the discussion here by sketching the Uhlenbeck compactification for $PU(2)$ monopoles, for details see \cite{FL2,T2} and \cite{Z1}. The important fact is that the quadratic map $K \circ \mu$ satisfies a pointwise properness condition,
\begin{equation*}
	\left( K(\mu(\Psi))\Psi, \Psi \right) \geq c^2 \,  \abs{\Psi}^4
\end{equation*}
for a real constant $c > 0$. This then yields, via the Weitzenb\"ock formula for the Dirac operator $D^+_A$, an `a-priori bound' for the spinor $\Psi$ of a monopole $[A,\Psi] \in M_{\kappa,\s}^{w}$, see for instance \cite{T3}:
\begin{equation} \label{a priori bound}
\norm{\Psi}^2_{\infty} \, \leq \, M/c^2 \ . 
\end{equation}
The positive constant $M$ only depends on fixed geometric data - the Riemannian metric, the fixed $Spin^c$ connection $B$, and the fixed automorphism $K$.

 This $C^0$ bound implies that there is an Uhlenbeck compactification of the moduli space $M_{\kappa,\s}^{w}$ by a space of `ideal monopoles', 
\begin{equation*}
	{\text{I}}M_{\kappa,\s}^{w} := \bigcup_{n=0}^{\infty} M^{w}_{\kappa-n,\s} \times \text{Sym}^n(X) \ ,
\end{equation*}
which is given a compact topology as in the instanton situation, and where each stratum has its previously defined topology. The closure $\overline{M}^{w}_{\kappa,\s} \subseteq \text{I}M_{\kappa,\s}^{w}$ is compact, see the corresponding statement in \cite{T2} for the analogous statement of \cite[Theorem 4.4.3]{DK}.

\subsection{Holonomy perturbations}
If we wish to obtain transversality for moduli spaces of $PU(2)$ monopoles we have to perturb the Dirac equation in the $PU(2)$ monopole equation, too. This is in contrast to the classical abelian Seiberg-Witten theory. There are different approaches in this situation \cite{FL,T2}, some using holonomy perturbations, some not.

 We wish to adopt the holonomy perturbations of the previous section to our situation. This is not so straight-forward to do with the Dirac equation if one wants to have an Uhlenbeck-compactification. In fact, we employ an artificial cut-off argument in order to keep a $C^0$ bound on the spinor component. However, we don't see how to control this cut-off term for a sequence of monopoles that tends to infinity. Even though the cut-off term converges for some subsequence, we cannot apply the same argument to get a $C^0$ bound on lower strata of the compactification. This is why we choose an additional cut-off factor to make the holonomy perturbation vanish for a sequence of monopoles converging to infinity. (Our equations in the lower strata will be perturbed in a `local' way without holonomy perturbations.)
\\

To a smooth connection $A \in \A$, a submersion $q: S^1 \times B \to X$ as in section \ref{Holonomy perturbations} and a complex-valued one-form $\alpha \in \Omega^1(X;\C)$ with support in $B$ we can associate a section
\begin{equation*}
 V_{q,\alpha}(A) := \alpha \tensor \Hol_q(A) \in \Omega^1(X;\gl(E)) \ .
\end{equation*}
There is then a result completely analogous to that of Proposition \ref{holonomy map}, the discussion following it, and the definition of the perturbation space $\W$: 

\begin{definition}
Let $\V$ be the Banach space of sequences ${\bf \alpha} = (\alpha_i)_{i \in \N}$, with each $\alpha_i$ an element of the Banach space $\Omega^1(B_i;\C)$, so that the sum 
\begin{equation*}
 \sum_i C_i \norm{\alpha_i}_{C^l(B_i)}
\end{equation*}
is finite. Here the constants $C_i$ are defined as in section \ref{Holonomy perturbations}.
For ${\bfalpha} \in \V$, the map $V_{\bfalpha}$ is defined to be the series $\sum V_{q_i,\alpha_i}$ which converges in the $C^\infty$ - topology of maps $\mathscr{A} \to L^2_{l}(X;\Lambda^1 \tensor \gl(E))$.
\end{definition}
The dependence of $V_{\bfalpha}$ on $\bfalpha$ is linear and we obtain a smooth map of Banach manifolds
\begin{equation*}
V: \V \times \A \to L^2_l(X;\Lambda^1 \tensor \gl(E)) \ , 
\end{equation*}
given by $({\bfalpha},A) \mapsto V_{\bfalpha}(A)$. This map is also equivariant with respect to the action of the gauge group $\G$ on $\A$ and on $\gl(E)$. 
\\

The perturbed $PU(2)$ monopole equations associated to the $Spin^c$ structure $\s$ and the Hermitian bundle $E \to X$, specified by the instanton number $\kappa$ and the line bundle $w$, and the perturbations $(\bfomega,\bfalpha) \in \V \times \W$ are then given by the equations

\begin{equation}\label{perturbed monopoles}
\begin{split}
   \dirac^+_{A} \Psi \, + \gamma(\beta) \Psi & = - h(\norm{F_A}_{L^m(X)}) h(\norm{V_{\bfalpha}(A)}_{L^\infty_{1,A}}) 
    \, \gamma(V_{\bfalpha}(A)) \Psi  \\
   \gamma(F_A^+) - K(\mu(\Psi)) \, & = -  \, h(\norm{F_A}_{L^m(X)}) \gamma(V_{\boldsymbol{\omega}}(A))  \ .
 \end{split}
\end{equation}
Here one has fixed some real number $m > 2$ and a positive function $h$ as required in Section \ref{moduli spaces with cut off perturbations}. 
%
The cut-off in the Dirac equation assures that the holonomy perturbation term $h(\norm{V_{\bfalpha}(A)}_{L^\infty_{1,A}}) 
    \, \gamma(V_{\bfalpha}(A))$ has uniformly bounded covariant derivative with respect to $A$. This yields a $C^0$ bound as in the classical case, using the Weitzenb\"ock formula and the maximum argument. For instance, the holonomy perturbations of Feehan-Leness naturally admit such bounds, and the proof of the $C^0$ bound in our situation is entirely analogous to theirs \cite[Lemma 4.4 and formula 2.23]{FL2}. Note also that the cut-off factors are gauge-invariant.

\begin{proposition}
Suppose $(A,\Psi)$ is a solution to the perturbed monopole equations (\ref{perturbed monopoles}). Then there is a positive constant $M'$, depending only on the Riemannian metric, the fixed $Spin^c$ connection $B$, the parameters $\beta$ and $K$ and the perturbation $(\bfomega,\bfalpha)$, so that there is a bound 
\begin{equation}\label{bound perturbed}
 	\norm{\Psi}_\infty \leq \frac{M'}{c^2} \ .
\end{equation}
\end{proposition}
The moduli space of $PU(2)$ monopoles perturbed by $(\bfomega,\bfalpha) \in \W \times \V$ is defined to be the space
\begin{equation*}
	M_{\kappa,\s}^{w}(\bfomega,\bfalpha):= \left\{ [A,\Psi] \in \bonf^{w}_{\kappa,\s}
	| (\ref{perturbed monopoles}) \text{ holds } \right\} \ .
\end{equation*}

Again, there is an elliptic deformation complex analogous to the one in Section \ref{properties of monopoles}. Let $(A,\Psi)$ be a solution to the perturbed monopole equations (\ref{perturbed monopoles}). Then we have the elliptic complex
\begin{equation} \label{elliptic complex}
	\begin{split}
	 0 \to L^2_{l+1}(X;\su(E)) & \xrightarrow{d^0_{(A,\Psi)}}  
L^2_l(X;S^+ \tensor E \oplus \Lambda^1 \tensor \, \su(E)) \\
& \xrightarrow{d^1_{(A,\Psi),(\bfomega,\bfalpha)}}
L^2_{l-1}(X;S^-\tensor E \oplus \Lambda^2_+ \tensor \, \su(E)) \to
0 \ .
	\end{split}
\end{equation}
The deformation operator is defined to be
\begin{equation}
 D_{(A,\Psi),(\bfomega,\bfalpha)} := d^{0,*}_{(A,\Psi)} \oplus d^1_{(A,\Psi),(\bfomega,\bfalpha)}\ .
\end{equation}
Note also that this deformation operator for the perturbed monopole equations and that for the unperturbed monopole equations only differ by addition of compact operators. In particular, their index is equal. As usual, the cohomology spaces of the above complex are denoted by $H^0_{(A,\Psi)}$, $H^1_{(A,\Psi),(\bfomega,\bfalpha)}$ and $H^2_{(A,\Psi),(\bfomega,\bfalpha)}$. 

\subsection{Compactification}\label{compactification of monopoles}

As in other considerations on $PU(2)$ monopoles, there is a natural Uhlenbeck-type compactification of the moduli space. As our holonomy perturbations are different from the ones used by Feehan-Leness \cite{FL2}, and as we have chosen a very specific cut-off behaviour of our holonomy perturbations towards infinity, there are a few essential modifications to be made here for obtaining the Uhlenbeck-compactification. We shall only give a sketch here by pointing out the essential differences from the arguments in \cite{FL2} and \cite{T}. All these approaches follow the ideas of \cite[Section 4.4]{DK}.\\

Recall that we have assumed $l\geq 3$. 
\begin{proposition}\label{compactness monopoles}
	Let $(A_n,\Psi_n)$ be a sequence of configurations, where $A_n$ are connections in the Hermitian bundle $E \to X$ and $\Psi_n$ are sections of the bundle $S^+ \tensor E$ , representing points $[A_n,\Psi_n]$ in the moduli space $M^{w}_{\s,\kappa}(\bfomega,\bfalpha)$. Then one of the following is true:
\begin{enumerate}
\item
After passing to a subsequence the sequence $[A_n,\Psi_n]$ converges in the $L^2_l$ topology of the moduli space $M^{w}_{\s,\kappa}(\bfomega,\bfalpha)$.\\
 	
\item
There is a point ${\bf x} \in {\text{\em Sym}}^n(X)$ for some $n\geq 1$, a connection $A'$ in a bundle $E' \to X$ and a section $\Psi'$ of the bundle $S^+ \tensor E'$, representing an element $[A',\Psi']$ of the moduli space $M^{w}_{\s, \kappa-n}(0,0)$ with the following property: After possibly passing to a subsequence, 
\\

\begin{enumerate}
\item[(i)]
there are bundle isomorphisms
\begin{equation*}
 	h_n: E|_{X \setminus {\bf x}} \to E'|_{X \setminus {\bf x}} \ ,
\end{equation*}
so that $(h_n)_*(A_n,\Psi_n)$ converges to $(A',\Psi')$ in $L^p_1(K)$ for all compact subset $K \subseteq X \setminus {\bf x}$, and for all $p \geq 2$, and 
\\

\item[(ii)]
the sequence of measures $(\abs{F_{A_n}}^2 vol_g)$ converges to the measure
\[ \abs{F_{A'}}^2 \text{vol}_g + 8 \pi^2 \sum_{x \in {\bf x}} \delta(x)\]
 in the weak-*-topology of measures, i.e. for any continous function $f$ on $X$ one has 
\[
	\int_X f\,  \abs{F_{A_n}}^2 \text{vol}_g \ \to \int_X f\, \abs{F_{A'}}^2 \text{vol}_g + 8 \pi^2 \sum_{x \in {\bf x}} f(x) \ 
\]
as $(n \to \infty)$.
\end{enumerate} 
\end{enumerate}
\end{proposition}

An ideal monopole of instanton number $\kappa$ is a pair $([A,\Psi],\bf{x})$ where ${\bf x} \in {\text{\em Sym}}^n(X)$, and $[A,\Psi]$ is an element of the unperturbed moduli space $M^w_{\s,\kappa-n}(0,0) = M^w_{\s,\kappa-n}$, if $n \geq 1$, and of the moduli space $M^w_{\s,\kappa}(\bfomega,\bfalpha)$ if $n=0$ . The space of $(\bfomega,\bfalpha)$-perturbed ideal monopoles of instanton number $\kappa$ is defined to be
\begin{equation*}
	{\text{I}}M_{\s,\kappa}^{w}(\bfomega,\bfalpha) := M_{\s,\kappa}^{w}(\bfomega,\bfalpha)\cup \bigcup_{n=1}^{\infty} M^{w}_{\s,\kappa-n} \times \text{Sym}^n(X) \ ,
\end{equation*}
with the notion of convergence between the different strata as in (2) of  Proposition \ref{compactness monopoles} (and the corresponding generalisation for sequences of ideal monopoles), and with each stratum having its original topology. It also follows from such kind of generalisation of this Proposition that the closure of the moduli space $M^{w}_{\s,\kappa}(\bfomega,\bfalpha)$ inside the space of ideal instantons ${\text{I}}M_{\s,\kappa}^{w}(\bfomega,\bfalpha)$ is compact. 
\\

In the following subsections we give a sketch of proof of this compactness result. We suppose that there is a sequence of configurations $(A_{n},\Psi_{n})$ representing points in the moduli space $M^{w}_{\s,\kappa}(\bfomega,\bfalpha)$. 

%



\subsubsection{Convergence of measures}
By the a priori bound (\ref{a priori bound}) the sequence $\Psi_{n}$ is bounded in $L^{\infty}$. By the second of the monopole equations we therefore see that the self-dual part of the curvatures $F^{+}_{A_{n}}$ 
is bounded in $L^{\infty}$. By the Chern-Weil formula there is therefore an $L^{2}$ bound on the sequence of curvatures $F_{A_{n}}$, and, in particular, the sequence of measures $\abs{F_{A_{n}}}^{2} \text{\em vol}_{g}$ is a bounded sequence of measures, considered as a sequence in the dual of $C^{0}(X;\R)$. The Banach-Alaouglu theorem claims that in the dual of a Banach-space the unit-ball is compact with respect to the weak-$^{*}$-topology, and this precisely the topology of the `convergence of measures'. Therefore there is a subsequence, without loss of generality equal to the sequence itself, that converges to a measure $\nu$ on $X$. There are only finitely many points $x \in X$ which do not lie in a geodesic ball of $\nu$-measure less than $\epsilon^{2}$, where $\epsilon$ will be the constant appearing in gauge fixing Uhlenbeck's theorem. In fact, there will be at most $\nu(X)/\epsilon^{2}$ many. Let ${\bf x} $ be the unordered collection of points where this occurs. 

\begin{remark}
If ${\bf x}$ is not the empty set we will have 
\begin{equation*}
	\lim_{n \to \infty} h(\norm{F_{A_{n}}}_{L^{m}(X)}) = 0
\end{equation*}
by Lemma \ref{L^m bound on curvature}, and because of the definition of the function $h$ in section \ref{moduli spaces with cut off perturbations}.
Therefore the holonomy perturbation terms in the perturbed monopole equations (\ref{perturbed monopoles}) tend to zero.
\end{remark}

\subsubsection{Compactness on the punctured manifold} We shall suppose now that ${\bf x}$ is not empty.
We begin by recalling Uhlenbeck's gauge fixing theorem \cite[Theorem 1.3]{U} for connections with sufficiently small curvature on the unit 4-ball with its standard metric.

\begin{lemma}\label{gauge fixing}
  There are constants $\epsilon, M \, > 0$ so that the following holds: Any connection $A$ of class $L^2_l$ on the trivial Hermitian rank-2 bundle over the closed unit ball $\overline{B}^4$ with $\norm{F_A}_{L^2} < \epsilon$ is gauge equivalent, by an element $u$ of class $L^2_{l+1}$, to a connection $u(A) = d + \tilde{A}$ so that the connection matrix $\tilde{A}$ satisfies
  \begin{equation*}
  \begin{split}
  		d^* \tilde{A} & = 0 \ \ , \\
		\norm{\tilde{A}}_{L^p_1} & \leq M \, \norm{F_{A}}_{L^p} \ \ 
  \end{split}
  \end{equation*}
for any $p \geq 2$.  
\end{lemma}
This is not precisely the formulation of the original source but it is rather standard to obtain the present formulation. Note also that by choosing $l \geq 3$ there is an inclusion $L^2_l \hookrightarrow L^p_1$ for any $p$. 

\begin{remark}
  The constants $\epsilon, M \, > 0$ in the preceding Lemma can be chosen so that the conclusion holds for any metric on the unit ball that is sufficiently close to the standard metric. We suppose from now on that such an $\epsilon$ is fixed. 
\end{remark}
On the open manifold $\Omega:= X \setminus {\bf x}$ any point $x \in \Omega$ admits a geodesic ball $B_x$ centered at $x$ so that 

\begin{equation}\label{bound curvature}
 \int_{B_x} \abs{F_{A_n}}^2 vol < \epsilon^2 
\end{equation}
for all $n \in \N$. 

\begin{lemma}\label{lemma L^p_1} Let $l\geq 3$. 
Provided the constant $\epsilon$ above is chosen small enough, any point $x \in \Omega$ admits a geodesic ball neighbourhood $B'_x$ centered at $x$ so that the following holds: 
\begin{enumerate}
\item
There is a sequence of gauge transformations $u_{n,x}$ of class $L^2_{l+1}$ over $B'_x$ so that the connections $u_{n,x}(A_n)$ are in Coulomb gauge with respect to some trivial connection over the ball $B'_x$,
\item the sequence of connections $u_{n,x}(A_n)$ is bounded in $L^p_1$ for any $p \geq 2$ by some constant that depends on $B'_x$ only, and
\item on the overlaps of the geodesic balls the sequence of gauge transformations $u_{n,x} u_{n,y}^{-1}$ is uniformly bounded in $L^p_2$ for any $p \geq 2$. 
\end{enumerate}
\end{lemma}
{\em Proof:}
After restricting to a possibly smaller ball $B'_x$ and rescaling the metric we may suppose that the metric is sufficiently close to the standard metric on the unit ball, and we may apply Uhlenbeck's gauge fixing theorem, Lemma \ref{gauge fixing}, as indicated above. 

Therefore there are gauge transformations $u_{n,x}$ so that the connection matrices $\tilde{A}_n$ of $u_{n,x}(A_n)$ with respect to some trivial connection are in Coulomb gauge, and there is a bound
\begin{equation} \label{estimate 1}
	\norm{\tilde{A_n}}_{L^2_1} \leq \, M \, \norm{F_{A_n}}_{L^2} \ ,
\end{equation}
with $M$ the constant of Lemma \ref{gauge fixing} above.
As the right hand side is uniformly bounded the sequence $\tilde{A_n}$ is uniformly bounded in $L^2_1$ norm. 

Over the ball $B'_x$ that we have identified with the standard 4-ball $B^{4}$ the second of the perturbed monopole equations may be rewritten, with slight abuse of notation, in the form
\begin{equation} \label{rewritten monopoles}
		d^+\tilde{A_n} = -(\tilde{A_n} \wedge \tilde{A_n})^+ 
			+ u_{n,x} \, \gamma^{-1} (K(\mu(\Psi_n))) \, u_{n,x}^{-1} 
			- u_{n,x} \, h(A_{n}) \, V_{\bfomega}(A_n) \, u_{n,x}^{-1} \ .
\end{equation}
By the a priori bound (\ref{bound perturbed}) and the construction of the holonomy perturbations, the second and third term of the right hand side of this equation are bounded in $L^{\infty}$. 

Now suppose we have $2 \leq p < 4$. Then we have bounded inclusions 
$L^p_1 \hookrightarrow L^q$ with $q = \frac{4p}{4-p}$ as well as a bounded multiplication map $L^4 \times L^q \to L^p$. We wish to apply elliptic regularity, so we transport the argument from the ball $B^{4}$ to the 4-sphere. Let $D \Subset B^{4}$ be an interior domain, and let $\psi : B^{4} \to [0,1]$ be a smooth cut-off function that is equal to $1$ on $D$ and that has compact support in the interiour of $B^{4}$. Let $\alpha_{n}:= \psi \, \tilde{A}_{n}$, that we may consider as a sequence of one-forms on the 4-sphere, containing $B^{4}$. Let us denote by $\chi_{n}$ the sequence given by the second and third term of (\ref{rewritten monopoles}). The one-forms $\alpha_{n}$ now satisfy the equation:

\begin{equation*}
\begin{split}
(d^{+}\oplus d^{*}) \alpha_{n} = (- (\alpha_{n} \wedge \tilde{A}_{n})^{+} + \psi \chi_{n},0) + (( d\psi \wedge \tilde{A}_{n})^{+}, *(d\psi \wedge * \tilde{A}_{n})) 
\end{split}
\end{equation*}
Therefore we get the following inequality by elliptic regularity of the operator $d^{+}\oplus d^{*}$ and the above mentioned bounded multiplication and inclusion maps
\begin{equation*}
\begin{split}
\norm{\alpha_{n}}_{L^{p}_{1}(S^{4})} & \leq c 		
	\left( \norm{\alpha_{n}}_{L^{p}_{1}(S^{4})} \norm{\tilde{A}_{n}}_{L^{4}} + \norm{\chi_{n}}_{L^{p}} + \norm{d\psi}_{L^{p}_{1}} \norm{\tilde{A}_{n}}_{L^{4}} \right)
\end{split}
\end{equation*}
for some positive constant $c$. Now using the estimate (\ref{estimate 1}) provided by Uhlenbeck's theorem we see that we can rearrange this inequality provided we choose $\epsilon$ in the bound (\ref{bound curvature}) small enough, so that we obtain a uniform bound on $\alpha_{n}$ in $L^{p}_{1}(S^{4})$ norm  for $2 \leq p < 4$. The desired unifom bound in $L^p_1(D)$ for $\tilde{A}_{n}$ now follows from 
\begin{equation*}
 \norm{\tilde{A}_{n}}_{L^{p}_{1}(D)} \leq \norm{\alpha_{n}}_{L^{p}_{1}(S^{4})} \ .
\end{equation*}
\\

We shall now extend the bound to arbitrary $p \geq 2$. In an intermediate step suppose $p \geq 3$. Then we certainly have $q > 8$. We obtain an inequality

\begin{equation*}
	\norm{\tilde{A}_n \wedge \tilde{A}_n}_{L^{q/2}} \leq \, \norm{\tilde{A}_n}_{L^q}^2 \leq c_p \norm{\tilde{A}_n}_{L^p_1}^2 \ , 
\end{equation*}
for some positive constant $c_p$. Now using the uniform $L^p_1$ bound on $\tilde{A}_n$ and repeating the above sequence of arguments we obtain (on a possibly smaller interiour domain) a uniform bound on $\tilde{A}_n$ in $L^{q/2}_1$ norms. Note, however, that for $r > 4$ there is a bounded inclusion $L^r_1 \hookrightarrow C^0 \subseteq L^\infty$, and as $q/2 > 4$ we are in this range now. Therefore we obtain bounds in arbitrary $L^p_1$ norms for arbitrary high $p \geq 2$ if we apply elliptic regularity and the equality (\ref{rewritten monopoles}) a third time. 
\\

The conclusion $(3)$ in the Lemma follows along the lines of `connections control gauge transformations' as in \cite[Section 2.3.7]{DK}. 

\qed

Next we shall `patch gauge transformations' to globalise results that hold on a countable union of balls $B_{x_i}$ to $\Omega$. The patching procedure of  gauge transformations is possible precisely because of (3) in the preceding Lemma, and the fact that for $p > 2$ there is a compact inclusion $L^p_2 \hookrightarrow C^0$, which allows us to pass to uniformly convergent subsequences from $L^p_2$ bounded sequences, and, therefore, we may make use of the exponential map for the gauge group as in \cite[Section 4.4.2]{DK}, \cite[Lemma A.0.1]{Fro2} and \cite[Section 3]{U} to globalise the gauge transformations. This yields the following 

\begin{proposition}
	In the situation above there is a subsequence $(n_{k})$, and there are gauge transformations $g_{n_{k}}$ on the punctured manifold so that the sequence $g_{n_k} (A_{n_k})$ is bounded in $L^p_1(K)$ for any compact subset $K \subseteq \Omega$, and any $p \geq 2$. 
\end{proposition}

Again, we shall assume that the subsequence is equal to the original sequence. 
\\

The construction of the holonomy perturbations involves a countable family of submersions $q_i:S^1 \times B_i \to X$. Now the points $y$ in the ball $B_i \subseteq X$ so that $q_i(z,y) \notin \Omega$ for some $z \in S^1$ form a subset of measure zero, and likewise the union of all such points in the family of balls $B_i$ forms a subset of measure zero. Therefore, the sections $V_{\bfomega}(g_n(A_n))$ and $V_{\bfalpha}(g_n(A_n))$ are well defined as $L^p(X)$ sections. The following Lemma appears as \cite[Lemma 10]{K}.

\begin{lemma}
There is a subsequence of $(n)$ for which the sequences $V_{\bfomega}(g_n(A_n))$, $V_{\bfalpha}(g_n(A_n))$ are Cauchy sequences as sections of class $L^p(X)$, for any $p \geq 2$. 
\end{lemma} 

From this fact and the perturbed Dirac equation, the first equation of (\ref{perturbed monopoles}), we can use elliptic regularity to show that the sequence $(g_n \Psi_n)$ is bounded in $L^p_1(K)$ for any compact subsets $K \subseteq \Omega$, and for all $p \geq 2$. So there is a  subsequence for which the spinor component $\Psi_{n}$ is Cauchy in $L^p(K)$ for all compact subset $K \subseteq \Omega$, and for all $p \geq 2$.  \\

From these global results we go back to a local consideration. Let us rename $(g_n(A_n),g_n \Psi_n)$ as $(A_n,\Psi_n)$. We restrict the monopole equations to the balls $B_{x_i}$ where the sequence of connections $A_n$ can be put in Coulomb gauge via a sequence of gauge transformations $v_{n,x_i}$. Then we use elliptic regularity to see that  the sequence $(v_{n,x_i}(A_n), v_{n,x_i} \Psi_n)$ is Cauchy in $L^p_1(B'_{x_i})$. Patching gauge transformations again,  there is a subsequence of $(n)$, without loss of generality equal to the original sequence, and there are gauge transformations $w_n$ on $\Omega$ so that $(w_n(A_n), w_n \Psi_n)$ is a Cauchy sequence in $L^p_1(K)$ for any compact subset $K \subseteq \Omega$, and any $p \geq 2$. Therefore, there is a connection $A$ on $E|_{\Omega}$ and a section $\Psi$ of $S^+\tensor E \, |_{\Omega}$ so that $(w_n(A_n),w_n \Psi_n)$ converges to $(A,\Psi)$ in any of the norms $L^p_1(K)$. We summarise our discussion in the following

\begin{proposition}
For the sequence $(A_n,\Psi_n)$ there is a subsequence, denoted without loss of generality by the original sequence, and there are gauge transformations $w_n$ on $\Omega$, of class $L^p_2(K)$ on each compact subset $K \subseteq \Omega$ and for each $p \geq 2$, so that the sequence $(w_n(A_n), w_n \Psi_n)$ converges to a configuration $(A,\Psi)$, defined on $\Omega$, in any of the norms $L^p_1(K)$. This configuration $(A,\Psi)$ solves the {\em unperturbed} monopole equations (\ref{pu2 monopoles}) over $\Omega$. 
\end{proposition} 

The next step consists of showing that  $(A,\Psi)$ extends to a configuration $(A',\Psi')$ of class $L^p_1$ on the entire manifold $X$, but associated to a different bundle $E' \to X$, and that this configuration solves the associated {\em unperturbed} monopole equations. 


\subsubsection{Regularity}

The regularity result that will be needed in the present situation is  the following:

\begin{prop}\label{regularity} (Regularity of $L^4$ small $L^2_1$ almost solutions) \\
There is a positive constant $M > 0$ so that the following holds: Let $(A,\Psi) \in L^2_1(Z;\su(E) \oplus S^+_\s \tensor E)$ be a $L^2_1$ configuration associated to a trivial Hermitian rank-2 bundle $E$ on the closed 4-manifold $Z$. Suppose that 
\begin{enumerate}
	\item $d^* A = 0$ \ , 
	\item $\norm{(A,\Psi)}_{L^4(Z)} \leq M \ \ $ , and 
	\item $(\dirac^+_A \Psi, \gamma(F_A^+) - \mu(\Psi))$ is of class $L^p$ for any $p \geq 2$. \\
\end{enumerate} 
Then $(A,\Psi)$ is of class $L^p_1$ for any $p \geq 2$. 
\end{prop}
{\em Proof:}
This follows from elliptic regularity along the same lines as in the Proof of Lemma \ref{lemma L^p_1}. For the first step one uses the same argument as in \cite[Proposition 3.2]{FL2} to obtain $L^p_1$ regularity for $2 \leq p < 4$ of $(A,\Psi)$, once its $L^4$ norm is small enough. 
\qed

\subsubsection{Removable singularities}
The following Lemma is analogous to the role of \cite[Proposition 4.4.10]{DK} in the compactification of the moduli space of ASD connections, see there also for the notion of `strongly simply-connectedness' of a domain of a manifold.

\begin{lemma} \label{strongly simply connected domain}
Let $\Omega$ be a strongly simply connected domain of a 4-manifold, and let $\Omega' \Subset \Omega $ be a precompact interior domain.  
There are constants $\delta_\Omega, \epsilon_\Omega > 0$ and $M_{\Omega,\Omega'} > 0$ so that the following holds: 

Let $A$ be a connection of class $L^p_1$ for some $p > 4$ on $\Omega$, satisfying 
\begin{enumerate}
	\item $\norm{F_A^+}_{L^p(\Omega)} < \delta_{\Omega}$ , and
	\item $\norm{F_A}_{L^2(\Omega)} < \epsilon_{\Omega}$ .
\end{enumerate}
Then $A$ can be represented over $\Omega'$ by a connection matrix $A^\tau$ 
with 
\begin{equation*}
	\norm{A^\tau}_{L^4(\Omega')} \leq M_{\Omega,\Omega'} \norm{F_A}_{L^2(\Omega)} \ .
\end{equation*}
\end{lemma}
{\em Proof:}
 The idea is similar to the proof of \cite[Proposition 4.4.10]{DK}. There is a cover of $\Omega'$ by a finite number of balls on which we can apply Uhlenbeck's gauge fixing theorem Lemma \ref{gauge fixing}, with possibly distinct values of $\epsilon$. The connection matrices $A^\tau, A^\sigma$ are related on their overlaps by transition functions $v$ satisfying
 \begin{equation*}
 	d v \, = \, v \, A^\tau + A^\sigma v \ .
 \end{equation*} 
As in Lemma \ref{lemma L^p_1} above we see that we can bound the $L^p_1$ norm of $A^\tau$ by its $L^4$ norm, and thus the $L^2$ norm of the curvature, and by the $L^p$ norm of $F_A^+$.  So by making these two quantities small we can make the $L^p_1$ norm of $A^\tau$ small. As we have a bounded inclusion $L^p_1 \hookrightarrow C^0$ we can therefore make the derivative of $dv$ small in $C^0$, by the preceding formula. One may then proceed as in the proof of \cite[Proposition 4.4.10]{DK} to obtain the conclusion.
\qed

\begin{remark}
Let $s > 0$ be a positive real number.
An easy computation shows that on an n-dimensional manifold $\Omega$ the $L^p$ norm of a differential form $\alpha$ behaves as in the following formula under  rescaling of the metric by the factor $s^2$:
\begin{equation*}
	\norm{\alpha}_{L^p(\Omega,s^2 g)} = s^{ \frac{n}{p} - \text{deg}(\alpha)} \norm{\alpha}_{L^p(\Omega,g)} \ .
\end{equation*}
Therefore, in the situation of the preceding Lemma, the $L^p$ norm of the (bundle-valued) 2-form $\chi$ for the metric $s^2 g$, with $s \in (0,1]$, controls its $L^p$ norm for the metric $g$, provided we have $p \geq 2$. As a consequence, if applying the previous Lemma with the role of $(\Omega,\Omega')$ played by annuli $(\mathscr{N},\mathscr{N'})$ inside the 4-ball, the constants in the Lemma may be taken independently of $r \in (0,1]$ under rescaling of the annuli by the dilatation map $x \mapsto r \, x$ of the ball.
\end{remark}
 
One now uses arguments analogous to the one in \cite[Section 4.4.5 and 4.4.6]{DK} for proving a Theorem on removable singularities, using Lemmata \ref{strongly simply connected domain} and \ref{regularity} at the place of Propositions \cite[4.4.10 and 4.4.13]{DK}. The remaining arguments leading to Proposition \ref{compactness monopoles} are analogous to \cite[Section 4.4]{T}. 


\subsection{Transversality in the main stratum}
We will now show that the perturbation space $\V \times \W$ is sufficient to obtain generic regularity (for the main stratum) of the moduli spaces $M^{w}_{\s,\kappa}(\bfomega,\bfalpha)$ of perturbed $PU(2)$ monopoles. Recall that the space $\conf^{**}$ consists of pairs $(A,\Psi)$ with irreducible connection $A$ and non-vanishing spinor $\Psi$. 
\begin{theorem}\label{transversality monopoles}
Suppose  the Condition \ref{condition} holds, and let $l \geq 3$. Then the smooth map of Banach manifolds 
\begin{equation*}
	\begin{split}
	\mathscr{F}: \V \times \W \times \conf^{**} \to L^2_{l-1}(X; S^- \tensor E \oplus \Lambda^2_+ \tensor \su(E)) \ ,
	\end{split}
\end{equation*}
given by the left-hand side of the perturbed $PU(2)$ monopole equations (\ref{perturbed monopoles}) minus the right hand side of (\ref{perturbed monopoles}), is transverse to zero. 
\end{theorem}
{\em Proof:} The proof is a generalisation to monopoles of the proof of the transversality theorem in \cite{K}.

Suppose we have an element $(\bfomega,\bfalpha,A,\Psi) \in \W \times \V \times \conf^{**}$ so that $\mathscr{F}((\bfomega,\bfalpha,A,\Psi))=0$. 
Let us denote by 
\begin{equation*}
\dirac^+_{A,\bfalpha} := \dirac^+_A + h(\norm{F_A}_{L^m(X)}) h(\norm{V_{\bfalpha}(A)}_{L^\infty_{1,A}}) \, \gamma(V_{\bfalpha}(A)) + \gamma(\beta)
\end{equation*}
 the `perturbed Dirac operator' in (\ref{perturbed monopoles}). We will show that the derivative
\begin{equation*}
\begin{split}
 	P:=d\mathscr{F}|_{(\bfomega,\bfalpha,A,\Psi)} : & \V \times \W \times L^2_{l}(X;S^+ \tensor E \oplus \Lambda^1 \tensor \, \su(E)) \\ & \to L^2_{l-1}(X; S^- \tensor E \oplus \Lambda^2_+ \tensor \su(E)) 
\end{split}
\end{equation*}
is surjective. To simplify notations we may assume that the `cut off factors' $ h(\norm{F_A}_{L^m(X)}) $ and $h(\norm{V_{\bfalpha}(A)}_{L^\infty_{1,A}})$ appearing in (\ref{perturbed monopoles}) are equal to $1$. It actually only matters that they are non-zero in this proof. In particular, we shall not consider the variation of these cut-off terms. Under this notational simplification, this derivative is given by the explicit expression
\begin{equation*}
	\begin{split}
	 	P(\bfnu,\bfbeta,a,\Phi)  = ( & \dirac^+_{A,\bfalpha} \Phi +
			\gamma( a +  \, dV_{\bfalpha}|_A (a) + \,  V_{\bfbeta}(A))\Psi 
			 , \\
			&  d_{A,\bfomega}^{+}a + \gamma^{-1} (K(\mu(\Psi,\Phi) + \mu(\Phi,\Psi))) + \, V_{\bfnu}(A) ) \ ,
	\end{split}
\end{equation*}
where $d^+_{A,\bfomega}$ now denotes the operator $d^+_A + \, dV_{\bfomega}|_{A}$.

Instead of proving directly that $P$ is surjective we shall consider the operators
\begin{equation*}
	P'_k: \W \times \V \times T_k \to L^2_{k-1} \ ,
\end{equation*}
where $T_k$ is the slice of the gauge-group action given by 
\[T_k=\ker(d_{A,\Psi}^{0,*}) \subseteq L^2_k(X;S^+ \tensor E \oplus \Lambda^1 \tensor \, \su(E)) \ , \] 
and where $P'_k$ has the same formal expression as $P$ above. 

As a first step we show that $P'_1:\W \times \V \times T_1 \to L^2$ is surjective. Suppose $(b,\Sigma) \in L^2(X;S^-\tensor E \oplus \Lambda^2_+ \tensor \, \su(E))$ is $L^2$-orthogonal to the image of $P'_1$. First we shall vary $\bfnu \in \W$ alone. We therfore have 
\begin{equation*}
 0 = \langle P'_1 (\bfnu) , (\Sigma,b) \rangle_{L^2} = \langle   V_{\bfnu}(A) , b \rangle_{L^2}
\end{equation*}
for all $\bfnu \in \W$. Now $\bfnu \mapsto V_{\bfnu}(A)$ has $L^2$-dense image (see the proof of the above Theorem \ref{transversality instantons} in \cite{K}) , so $b=0$. 

Next we would like to vary the spinor $\Phi$ alone, i.e. to consider $P'_1(\Phi)$ for an arbitrary $L^2_1$ section in $S^+ \tensor E$. This, however, is not possible because only the slice $T_1$ is involved in the definition of $P'_1$ -- it is not clear whether for a general spinor $\Phi$ there is a solution $(0,\Phi) \in T_1$. If it were possible, the argument would continue like this: We would have
\begin{equation*}
 0 = \langle P'_1 (\Phi) , (\Sigma,0) \rangle_{L^2} = \langle \dirac^+_{A,\bfalpha} \Phi, \Sigma \rangle_{L^2} 
\end{equation*}
for all $\Phi \in L^2_{1}(X;S^+ \tensor E)$. As a consequence, we would get $\dirac^-_{A,\bfalpha} \Sigma = 0$ in the distributional sense, where $\dirac^-_{A,\bfalpha}$ is the formal $L^2$-adjoint of $\dirac^+_{A,\bfalpha}$. By elliptic regularity of the operator $\dirac^-_{A}$ we would then see that $\Sigma$ is actually of Sobolev class $L^2_{l}$, so $\dirac^-_{A,\bfalpha} \Sigma = 0$ would hold in the usual sense. Now the point is that we still can conclude that $\Sigma$ satisfies the Dirac equation $\dirac^-_{A,\bfalpha} \Sigma = 0$. Here is the argument: First, there is an elliptic deformation complex as in (\ref{elliptic complex}) for any Sobolev index $k$, in particular for $k=1$. Applying Hodge theory to this elliptic complex gives the topological decomposition:
\begin{equation*}
 L^2_{1}(X;S^+ \tensor E \oplus \, \Lambda^1 \tensor \, \su(E)) 
	= \ker(d^{0,*}_{(A,\Psi)} \oplus \text{im}(d^0_{(A,\Psi)}) = T_1 \oplus \text{im}(d^0_	{(A,\Psi)}). 
\end{equation*}
Second, we observe that $\im(d^1_{(A,\Psi),(\bfomega,\bfalpha)}) = \im(d^1_{(A,\Psi),(\bfomega,\bfalpha)} |_{T_1})$ because $d^1_{(A,\Psi),(\bfomega,\bfalpha)} \circ d^0_{(A,\Psi)} = 0$. Third, the restriction of $P'_1$ to the slice $T_1$ is precisely equal to $d^1_{\bfomega,\bfalpha} |_{T_1}$. As a consequence, for any $\Phi \in L^2_{1}(X;S^+ \tensor E)$ there is an element $(a,\Phi') \in T_1$ of the slice so that  $P'_1(a,\Phi') = (\dirac_{A,\bfalpha}^+ \Phi, \dots)$ (we are only interested in the spinor component of $P'_1$ in the argument to show that $\Sigma$ satisfies the Dirac equation).

By assumption $\Psi \neq 0$, so if $\Psi(x) \neq 0$ then $\Psi(y)$ is non-zero for all $y$ in a neighbourhood $U$ of $x$. We shall now vary $\bfbeta$ alone, so that we obtain
\begin{equation}\label{vanishing of spinor}
 0 = \langle P'_1(\bfbeta) , (\Sigma,0) \rangle_{L^2} = \langle  \gamma(V_{\bfbeta}(A)) \Psi , \Sigma \rangle_{L^2}
\end{equation}
for all $\bfbeta \in \V$. Note that the map $\Lambda^1 \tensor \gl(E) \to S^- \tensor E$ given by $ e \mapsto \gamma(e) \Psi$ is pointwise surjective at any point where $\Psi \neq 0$. By the condition \ref{condition} and the fact that $A$ is irreducible there is a finite number of submersions $q_i:S^1 \times B_i \to X$ so that $x \in \text{int}(B_i)$ and so that the holonomys $\Hol_{q_i}(A) (x) \in \gl(E)_x$ span $\gl(E)_x$. Furthermore, the holonomy sections $\Hol_{q_i}(A)$ on $B_i$ continue to span $\gl(E)$ in a neighbourhood of $x$, because $l$ was chosen so that we have a Sobolev inclusion $L^2_l \hookrightarrow C^0$, and so the holonmy sections are continous. Multiplying the finite number of holonomy sections $\Hol_{q_i}(A)$ with convenient one-forms $\beta_i \in \Omega^1(B_i;\C)$ supported in small enough neighbourhoods of $x$ we obtain a perturbation $\bfbeta=(\beta_i) \in \V$, where all but these finite number of one-forms are zero, so that the equation (\ref{vanishing of spinor}) implies that $\Sigma$ is zero in a neighbourhood $U'$ of $x$. As we also have $\dirac^-_{A,\bfalpha} \Sigma = 0$ the unique continuation principle \cite{A} for solutions to the perturbed Dirac-equation implies that $\Sigma = 0$ on the whole of $X$. 

Therefore we have shown that $P'_1$ is surjective. Suppose now that $P'_1(\bfbeta,\bfnu,b,\Phi)$ lies in the Sobolev class $L^2_{l-1}$. The additional hypothesis that $(b,\Phi) \in T_1$ now imply by elliptic regularity that $(b,\Phi)$ is of Sobolev class $L^2_{l}$. Obviously we then have $P'_l(\bfbeta,\bfnu,b,\Phi)= P'_1(\bfbeta,\bfnu,b,\Phi)$, so that $P'_l$ and in particular $P$ is in fact surjective onto $L^2_{l-1}$. \qed

\begin{corollary}\label{generic regularity monopoles}
	For a residual set of perturbations $(\bfomega,\bfalpha) \in \V \times \W$ the subspace $M^{w,**}_{\s,\kappa}(\bfomega,\bfalpha)$ is regular. It therefore admits the structure of a smooth manifold of the expected dimension. 
\end{corollary} $\hfill \square$

\begin{remark}
This transversality Theorem remains true with literally the identical proof if the bundle $E \to X$ is replaced with any Hermitian bundle of higher rank than 2, contrary to the situation in \cite{FL2}. This may be relevant to any investigation following the ideas developed by the author in  \cite{Z1}.  
\end{remark}

\begin{lemma} \label{surjectivity dirac} Suppose the connection $A$ is irreducible and that the Dirac-operator $\dirac^+_A$ has non-negative index. Then there is an open and dense subset of elements $\bfalpha \in \V$ so that the perturbed Dirac-operator
\begin{equation*}
\begin{split}
 \dirac^+_{A,\bfalpha} & = \dirac^+_A \, + \, h(\norm{F_A}_{L^m(X)}) h(\norm{V_{\bfalpha}(A)}_{L^\infty_{1,A}}) \, \gamma(V_{\bfalpha}(A)) + \gamma(\beta) \\
 &  : L^2_{l}(X;S^+ \tensor E) \to L^2_{l-1}(X;S^- \tensor E)
\end{split}
\end{equation*}
is surjective.
\end{lemma}
{\em Proof:}
Let us consider the map 
\begin{equation*}
\begin{split}
  g: \V \times L^2_l(X;S^+\tensor E) \setminus \{ 0 \} & \to L^2_{l-1}(X;S^- \tensor E) \\
	({\bfalpha},\Psi) &\mapsto \dirac^+_{A,\bfalpha} \Psi \ .
\end{split}
\end{equation*}
As in the proof of the last theorem we see that $0$ is a regular value of this map. Let $\mathscr{M}$ be the zero-set $g^{-1}(0) \subseteq \V \times L^2_l(X;S^+\tensor E) \setminus \{ 0 \}$. The projection onto the first factor $\pi: \mathscr{M} \to \V$ is then a Fredholm map of the same index as that of the Dirac operator $\dirac^+_A$. In fact, for $(\bfalpha,\Psi) \in \mathscr{M}$ the kernel and cokernel of $d \pi_{(\bfalpha,\Psi)}$ and of $\frac{\partial g}{\partial \Psi}|_{(\bfalpha,\Psi)}$ are naturally isomorphic. Now by the Sard-Smale theorem there is a residual subset of $\V$ consisting of regular values for $\pi$. Note that we simply have 
\begin{equation*}
 \left.\frac{\partial g}{\partial \Psi}\right|_{(\bfalpha,\Psi)} = \dirac^+_{A,\bfalpha} \ .
\end{equation*}
Therefore, if $\bfalpha$ is a regular value for $\pi$, the perturbed Dirac-operator $\dirac^+_{A,\bfalpha}$ is surjective. The dependence of the bounded operator $\dirac^+_{A,\bfalpha}: L^2_{l}(X;S^+ \tensor E) \to L^2_{l-1}(X;S^- \tensor E)$ on $\bfalpha$ is continous. Therefore the residual set of values $\bfalpha$ for which this operator is surjective is also open.  \qed

\subsection{Transversality in the lower strata}
Recall that the unperturbed monopole equations (\ref{pu2 monopoles}) involves the `parameters'
$\beta \in \Omega^1(X;\C)$ and $K \in \, \text{Gl}_{+}\, (\Lambda^2_+)$. We did not refer to these as `perturbations' of the equations in order to reserve this term for the holonomy perturbations $(\bfomega,\bfalpha)$ considered above (and to avoid confusion). Nonetheless we use these parameters to obtain regularity for the (holonomy) unperturbed moduli spaces $M^{w,**}_{\s,\kappa-n}$ for $n \geq 1$, following \cite{T}. A similar approach appears in \cite{F}. 

In fact, the whole point about the holonomy perturbations was to extend the perturbations appearing in (\ref{perturbed asd cut off}) which defined the Casson-type invariant $n_{w,o}(\bfomega)$, and to obtain Lemma \ref{surjectivity dirac} above. 

Regularity in \cite{T} is achieved by first perturbing the equations with $(\beta,K)$, resulting in a `parametrised moduli space' which contains exceptions to transversality, and then perturbing the Clifford-multiplication and thereby the Riemannian metric to remove these points of exception. The equations (\ref{pu2 monopoles}) are therefore considered as equations for the configuration $(A,\Psi)$ and the parameters $(\beta,K,\gamma)$. After fixing a reference Clifford multiplication, $\gamma$ can be seen as an orientation preserving automorphism of the cotangent bundle $\Lambda^{1}(X)$. Correspondingly, the perturbation space $\mathscr{P}^{l}$ is identified with the space
\[
	\Omega^{1}(X;\C) \times \Gamma( X; \text{Gl}_{+} \, (\su(S^{+}_{\s})) \times \Gamma( X; \text{Gl}_{+} \, (\Lambda^{1}(X)) \ ,
\]
completed with respect to the $C^{l}$ norm. 

The regularity theorem is then the following, where the term `abelian locus' refers to configurations $(\Psi,A)$, where the connection $A$ yields a parallel splitting $E = L_1 \oplus L_2$ into two line bundles, and where the {\it non-vanishing} spinor $\Psi$ is a section that lies entirely in one of the resulting $U(2)$ bundles $S^+_\mathfrak{s} \otimes L_i$. This is best understood in the situation of the $S^1$ action described in Section \ref{circle action} below.

\begin{theorem}\cite[Theorem 3.19]{T}\label{regularity lower strata}
There is a dense second category subset $\mathscr{P}^{l}_{0}$ of $\mathscr{P}^{l}$ so that for any $(\beta,K,\gamma) \in \mathscr{P}^{l}_{0}$ the corresponding moduli space $M^{w,**}_{\s,\kappa} = M^{w,**}_{\s,\kappa}(\beta,K,\gamma)$ is regular away from the abelian locus. Moreover, the space $ \mathscr{P}^{l}_{0}$ may be chosen so that this conclusion holds for all moduli spaces $M^{w,**}_{\s,\kappa-n}$ with $n \geq 0$ at once.
\end{theorem}

\subsection{Index computations}
The index of the elliptic
operator $\dirac_A^+$ can be computed from the Atiyah-Singer index theorem, and is given by
\begin{equation*}
 \text{ind}_\C(\dirac_A^+) = \langle \text{ch}(E) e^{\frac{1}{2} c_1(S^+_\mathfrak{s})} \hat{\text{A}}(TX) ,
[X]
\rangle \ .
\end{equation*}
This formula can be found in \cite[\textsection6.4]{S}. In our situation we thus obtain
\begin{equation} \label{index, explicit}
 \text{ind}_\C(\dirac_A^+) = 2 \, \text{ind}_\C (\dirac^+_\mathfrak{s}) + \frac{1}{2}  \, \langle
c_1(S^+_\mathfrak{s}) c_1(E) + c_1(E)^2 ,[X] \rangle - \langle c_2(E),[X] \rangle \ ,
\end{equation}
where we have denoted by $\dirac^+_\s$ the Dirac operator determined by the (untwisted) $Spin^c$ structure $\s$ and the
fixed $Spin^c$ connection that we have suppressed in the notation. 

If we wish to make explicit that the coupled Dirac operator $\dirac_A^+$ depends on the $Spin^c$ structure
$\s$ and the Hermitian bundle $E \to X$, then we shall write $\ind(\dirac_A^+)=: \ind(\dirac_\s^+,E)$.
Similarly, we shall write $\ind(\delta_A)=:\ind(\delta,E)$ for the deformation operator $\delta_A$ of the
instanton moduli space. \\

For our application in mind the bundle $E \to X$ is chosen so that it defines a Casson-type instanton
moduli space as in Definition \ref{casson invariant perturbation}. This is only a restriction up to tensoring the
Hermitian bundle $E \to X$ with a Hermitian line bundle, see Proposition \ref{tensor product with line
bundle}. We are, however, free in the choice of the $Spin^c$ structure.

%

\begin{remark}
For the (untwisted) Dirac operator associated to a $Spin^c$ structure $\s$ on a negative definite four-manifold, the index is always non-positive:
\begin{equation*}
	\ind_\C (\dirac^+_{\s}) = \frac{1}{8} \langle c_1(S^+_\s)^2 ,[X] \rangle - \text{\em sign}(X) \leq 0 \ .
\end{equation*}
\end{remark}
The following Proposition shows that this is not necessarily the case for the twisted Dirac operator.
\begin{prop} \label{index} Let $X$ be a negative definite four-manifold, and suppose $E \to X$ is a Hermitian rank-2 bundle. Suppose one has 
\begin{equation}\label{condition kappa}
	0 = \langle c_2(E) - \frac{1}{4} c_1(E)^2, [X] \rangle  \ .
\end{equation}
Then one can find a $Spin^c$ structure $\s$ so that the index of the twisted Dirac operator $\dirac^+_A :
\Gamma(X;S^+_\s \tensor E) \to \Gamma(X;S^-_\s \tensor E)$ is one,
\begin{equation*}
 \ind_\C(\dirac^+_\s,E) = 1 \ .
\end{equation*}
In particular, in the situation of a bundle $E$ leading to our Casson-type instanton invariant defined in Definition \ref{casson invariant perturbation}, this can always be achieved.
\end{prop}
{\em Proof:}
Notice first that if we take the tensor product of $E \to X$  with a Hermitian line bundle $L \to X$, then
we have 
\begin{equation*}
 \ind (\dirac^+_{\s \tensor L^*},E \tensor L) \, = \, \ind (\dirac^+_{\s},E) \ ,
\end{equation*}
where $\s \tensor L^*$ denotes the $Spin^c$ structure that one obtains from twisting $\s$ with the dual of
the line bundle $L$. Note also that tensoring $E$ with a line bundle leaves the condition (\ref{condition kappa}) unchanged. 

Suppose we are given a set of elements $\{ e_i \}$ inducing a basis of $H^2(X;\Z)/\text{\em
Torsion}$ which diagonalises the intersection form. Up to taking tensor products with line bundles, and up
to a torsion element, we may suppose that 
\begin{equation*}
 c_1(E) = \sum k_i \, e_i  \ , 
\end{equation*}
with $k_i$ being either 0 or 1. In fact,  for a line bundle $L$ we have $c_1(E \tensor L) = c_1(E) + 2 c_1(L)$.  Up to a permutation of the indices, we may therefore suppose that $k_1= \dots = k_N = 1$, and
$k_{N+1} = \dots = k_{b_2(X)} = 0$.

By our liberty of choice in the $Spin^c$ structure $\s$ we may suppose
that $c_1(\s)$ is any characterstic element of the intersection form. In particular, we must have
\begin{equation*}
 c_1(\s) = \sum l_i \, e_i 
\end{equation*}
with $l_i \equiv 1 \ (\text{mod} \ 2)$. Let us take $\s_0$ to be a $Spin^c$ structure so that $c_1(\s_0) = \sum_{i=1}^{b_2(X)} e_i$. One easily checks that 
\begin{equation*}
\ind_\C(\dirac^+_{\s_0},E)  = -N - \langle \,  c_2(E), [X] \, \rangle = -\frac{3N}{4} \ . 
\end{equation*}

We may get any other $Spin^c$ structure $\s$ by tensoring $\s_0$ with  a line bundle $K$, and we have $c_1(\s_0 \tensor K) = c_1(\s_0) +2 c_1(K)$. 

Suppose we have $c_1(K) = \sum_{i=1}^{b_2(X)} m_i e_i. $ Then the formula (\ref{index, explicit}) implies that 
\begin{equation}\label{index, detailed}
	\begin{split}
	\ind_\C(\dirac^+_{\s_0 \tensor K},E) & = \ind_\C(\dirac^+_{\s_0},E) \\
			&  \ \ + \langle c_1(\s_0) c_1(K) + c_1(K)^2 + c_1(K) c_1(E), [X] \, \rangle \  \\ 
			& = -\frac{3N}{4} - \sum_{i=1}^{N} (2m_i + m_i^2) - \sum_{i=N+1}^{b_2(X)} (m_i + m_i^2) \ . 
	\end{split}
\end{equation}
Now $2m_i + m_i^2 \geq -1$, and for $m_i= -1$ we have $2m_i +m_i^2 = -1$. Let $l:= 1 +3N/4$. Then we see from this formula that if we set $m_i = -1$ for $i=1, \dots, l$ and $m_i=0$ for $i= l+1, \dots, b_2(X)$, we do indeed obtain $\ind_\C(\dirac^+_{\s_0 \tensor K},E)=1$. 

\qed 
\begin{remark}
	The conclusion of the preceding Proposition certainly may hold under weaker assumptions on the intersection form of $X$ and the bundle $E \to X$, respectively the number $\langle p_1(\su(E)), [X] \rangle$. This puts the techniques developed in this paper in the perspective of further applications. However, we shall content ourselves here with the version above which is sufficient for our needs. 
\end{remark}

\begin{proposition}\label{expected dimensions compactification}
Suppose $X$ is a negative definite four-manifold with $b_1(X) = 1$, and $E \to X$ is a bundle leading to the Casson-type instanton moduli space defined in Definition \ref{casson invariant perturbation}.

	Suppose a $Spin^c$ structure $\mathfrak{s}$ is chosen so that the index of the twisted Dirac operator is one, which is possible according to the preceding Proposition. Then the expected dimension of the moduli
space of $PU(2)$ monopoles $M_{0,\s}^w$ is 2. Furthermore, the lower strata of the space of ideal monopoles
$\text{I}M_{0,\s}^{w} = \bigcup_{n=0}^{\infty} M_{-n,\s}^{w}$ are all of strictly negative expected
dimension. 
\end{proposition}
{\em Proof:} From the formulae (\ref{ex-dim instantons}) and (\ref{index, explicit}) the expected dimension
of the moduli space
$M_{-n,\s}^{w}$ equals 
\begin{equation*}
\begin{split}
 \ind(\delta,E_{-n}) + 2 \, \ind(\dirac^+_\s,E_{-n}) & = \ind(\delta,E) - 8n +  2\, (\ind(\dirac^+_{\s},E) +
n) \\
 		& = \ind(\delta,E) + 2 \, \ind(\dirac^+_\s,E) - 6n .\\
\end{split}
\end{equation*}
For $n=0$ the assumptions made imply that this number is equal to 2. For $n > 0$ this number is $2-6n$, and so strictly negative. \qed

\subsection{A circle action on the configuration space of $PU(2)$ monopoles}\label{circle action}
We will now introduce an action of the circle $S^1$ on the configuration space $\bonf^w_{\kappa,\s}$ of
$PU(2)$ monopoles, and we will describe its fixed point set. On the pre-configuration space this action $S^1
\times \conf \to \conf$ is simply given by scalar multiplication on the spinor,
\begin{equation*}
	\left( z, (A,\Psi)\right) \mapsto (A,z \Psi) \, .
\end{equation*}
This action descends to an action $\rho: S^1 \times \bonf^{w}_{\kappa,\s} \to \bonf^{w}_{\kappa,\s}$ on the
configuration space. However, the latter is not effective -- in fact it is the two-fold covering of an
effective action $\rho^{1/2}$ because the stabiliser $\Gamma_A$ of a connection $A$ in the
gauge groupe $\G$ always contains the centre $Z(SU(2)) = \pm \id$. We therefore have $\rho(z,[A,\Psi]) = [A,
z \Psi]$ and $\rho^{1/2}(z,[A,\Psi]) = [A,z^{1/2} \Psi]$, where $z^{1/2}$ is an arbitrary square-root of
$z$. The following well-known proposition can be found, for instance, in \cite{FL3,T,Z1}. 

\begin{proposition}\label{fixed points}
A configuration $[A,\Psi] \in \bonf^{w}_{\kappa,\s}$ belongs to the fixed point set of $\rho$ respectively
$\rho^{1/2}$ if and only if for any representative $(A,\Psi) \in \conf$ we have one of the following:
\begin{enumerate}
\item There is an A-parallel decomposition of $E$ into the sum of two line-bundles, $E = K \oplus L$, and
$\Psi$ is a non-vanishing section of either $S^+ \tensor K$ or $S^+ \tensor L$. 
\item The spinor vanishes, $\Psi \equiv 0$.
\end{enumerate}
On the complement of the fixed point set the action $\rho^{1/2}$ is free. 
Furthermore, if $(A,\Psi)$ solves the monopole equations (\ref{pu2 monopoles}) then $[A,\Psi]$ is a fixed point if and only if $A$ is reducible or has vanishing spinor component. In particular, the action is free on the subspace $M^{w,**}_{\s,\kappa}$ of $M^{w}_{\s,\kappa}$.
\end{proposition}
\begin{definition}\label{abelian locus}
	Configurations $[A,\Psi]$ in the fixed point set with {\bf non-vanishing spinor component $\Psi$}, necessarily of type $(1)$ above, are referred to as the `abelian locus'.
\end{definition}

\begin{remark}
On negative four-manifolds with the bundle $E \to X$ chosen as in Lemma \ref{no reductions} there are no
fixed points of the first type in the above Proposition. This applies to the situation in which the Casson-type invariant is defined, see Definition \ref{casson invariant perturbation}.
\end{remark}

\subsection{Orientations}\label{orientations monopoles}
In analogy to the instanton situation we will determine natural orientations for the subspace of the $PU(2)$ moduli space $M_{\kappa,\s}^{w}(\bfomega,\bfalpha)$ consisting of regular points with {\em trivial} stabiliser under action of the gauge group, and thus forming a manifold of the expected dimension. Again, this is done by considering the determinant line bundle of the family of Fredholm operators given by the deformation operators parametrised by the configuration space $\bonf^{w,*}_{\kappa,\s}$. 

Let 
\[D_{(A,\Psi),\bfomega,\bfalpha} :=d^{0,*}_{(A,\Psi)} \oplus d^1_{(A,\Psi),(\bfomega,\bfalpha)}
\] 
be the deformation operator associated to the `elliptic complex' (\ref{elliptic complex}). To be more precise, (\ref{elliptic complex}) is only a complex if $(A,\Psi)$ satisfies the perturbed monopole equations. However, the operator $D_{(A,\Psi),\bfomega,\bfalpha}$ is a real Fredholm operator for any $(A,\Psi) \in \conf$. We obtain a homotopy of Fredholm operators by the formula $t \mapsto D_{(A,t \Psi),(\bfomega,\bfalpha)}$. For $t=0$ this operator takes the form of a direct sum:
\begin{equation}\label{direct sum fredholm}
D_{(A,0),\bfomega,\bfalpha} = \delta_{A,\bfomega} \oplus \dirac_{A,\bfalpha}^+
\end{equation}
Let
\begin{equation*}
	\Theta_{\bfomega,\bfalpha}(t) \to \bonf^{w,*}_{\kappa,\s}
\end{equation*}
be the quotient of the determinant line bundle of the Fredholm operator $D_{(A,t\Psi),\bfomega,\bfalpha}$ by the gauge group $\G$. Likewise, let 
\begin{equation*}
\Omega_{\bfalpha} \to \bonf^{w,*}_{\kappa,\s}
\end{equation*}
be the quotient of the real determinant line bundle of the family of Dirac operators $\dirac_{A,\bfalpha}^+$, considered as $\R$-linear Fredholm operators,  by the gauge group $\G$.

We can now state the following Lemma:
\begin{lemma}
\begin{enumerate} 
\item	Suppose the moduli space $M^{w,**}_{\kappa,\s}(\bfomega,\bfalpha) \subseteq \bonf^{w,**}_{\kappa,\s}$ is regular. Then the restriction of the determinant line bundle $\Theta_{\bfomega,\bfalpha}(1)$ to it is equal to the orientation bundle of $M^{w,**}_{\kappa,\s}(\bfomega,\bfalpha)$. 
\item The restriction of the determinant line bundle $\Theta_{\bfomega,\bfalpha}(0)$ to the subspace $\bonf^{w,**}_{\kappa,\s}$ is given by
\begin{equation*}
	\Theta_{\bfomega,\bfalpha}(0)|_{\bonf^{**}} = \pi^* \Lambda_{\bfomega} \tensor \Omega_{\bfalpha} \ ,
\end{equation*}
where $\Lambda_{\bfomega} \to \bonf^{w,*}_{\kappa}$ is the determinant line bundle of the instanton deformation operator of section \ref{orientations instantons}, the map $\pi : \bonf^{w,**}_{\kappa,\s} \to \bonf^{w,*}_{\kappa}$ is the projection onto the connection component, and where $\Omega_{\bfalpha}$ is the real determinant line bundle of the family of Dirac operators $\dirac_{A,\bfalpha}^+$, considered as $\R$-linear operators.  
\end{enumerate}
\end{lemma}
{\em Proof:}
The first statement follows directly from the definition. The second statement is clearly a consequence of formula (\ref{direct sum fredholm}) above. 
\qed

Now note that the line bundles $\Theta_{\bfomega,\bfalpha}(1)$ and $\Theta_{\bfomega,\bfalpha}(0)$ are naturally isomorphic, at least up to multiplication by a nowhere vanishing positive real function. In fact, they occur as restrictions of a line bundle over $[0,1] \times  \bonf^{w,*}_{\kappa,\s}$ to $\{ 0 \} \times \bonf^{w,*}_{\kappa,\s}$ respectively $\{ 1 \} \times \bonf^{w,*}_{\kappa,\s}$. So, up to multiplication by a non-vanishing positive function, there is a unique non-vanishing section of the bundle over each stripe $[0,1] \times \{ [A,\Psi] \}$ for each $[A,\Psi] \in \bonf^{w,*}_{\kappa,\s}$. 

Note further that the line bundle $\Omega_{\bfalpha}$ is trivial. To see this, notice that there is a natural inclusion $Gl(n,\C) \hookrightarrow Gl(2n,\R)$ compatible with the convention, that if $(v_1, \dots, v_n)$ is a basis of $\C^n$, then $(v_1, i v_1, \dots, v_n, iv_n)$ is the naturally associated basis of $\C^n$ considered as $\R^{2n}$. The complex entries of an $n \times n$ matrix then become blocks of $2 \times 2$ matrices with real entries. It is well-known that this inclusion $Gl(n,\C) \hookrightarrow Gl(2n,\R)$ factors through $Gl_+(2n,\R)$, the subgroup of $Gl(2n,\R)$ with positive determinant -- this is what is meant with `a complex vector space is canonically oriented when considered as the underlying real vector space'. As a consequence, the real determinant line bundle of a complex vector bundle, seen as real vector bundle, is trivial. The given convention here also determines a natural trivialisation of the real determiant line bundle: If, in a fibre, a complex basis is given by $(v_1, \dots, v_n)$, then $v_1 \wedge i v_1 \wedge \dots \wedge v_n \wedge i v_n$ is a basis of the real determinant line, and this definition is independent, up to a positive multiple, of the chosen basis.

Now, the  Fredholm operators $\dirac_{A,\bfalpha}^+$ are complex linear operators, and $\G$ acts complex linearly on the spinor component. The preceding discussion applies and one concludes that $\Omega_{\bfalpha}$ is indeed trivial with natural trivialisation.


Recall that the letter $o$ designated a choice of trivialisation of the determinant line bundle $\Lambda_0 \to \bonf^{w,*}_{\kappa}$ in section \ref{orientations instantons}, and that there is also a natural isomorphism between $\Lambda_0$ and $\Lambda_{\bfomega}$. We therefore have proved the following

\begin{corollary}\label{natural orientation monopoles}
	Suppose the moduli space $M^{w,**}_{\kappa,\s}(\bfomega,\bfalpha)$ is regular. 
	Then it is orientable. Furthermore, a choice of trivialisation $o$ of the trivial line bundle $\Lambda_0$ determines a natural orientation of $M^{w,**}_{\kappa,\s}(\bfomega,\bfalpha)$.
\end{corollary} 
\qed

\subsection{Local models around the instantons}\label{local model}
We will recall some theory of local models in general here, and then apply the results to the neighbourhood of the instantons $[A] \in M_{0}^{w}(\bfomega)$ inside the $PU(2)$ monopole moduli space $M_{0,\s}^{w}(\bfomega,\bfalpha)$, where the instantons are considered as $PU(2)$ monopoles $[A,0]$ with vanishing spinor. We actually need some equivariant version for local models here, and we will emphasise why this works as well. Furthermore, we shall make use of the elliptic deformation complex (\ref{elliptic complex}). To simplify the notations we shall continue to write $H^0, H^1, H^2$ if we mean actually the {\em harmonic representatives} of the cohomology spaces $H^0, H^1$ and $H^2$ of that complex. 
\\

The following results on local models are standard and are usually refered to as `Kuranishi models', see for instance \cite[Section 4.2.5]{DK} or the corresponding description for $PU(2)$ monopoles in \cite{FL2}. To set up our notation, we shall make this explicit in the equivariant setting, even though this equivariant model can be found, for instance, in Freed and Uhlenbeck's treatise as \cite[Lemma 4.7]{FU}.

\begin{proposition}\label{slices}
Suppose $(A,\Psi)$ is an element of the pre-configuration space $\conf$. Let $T_{(A,\Psi)}$ be the slice of the gauge-group action $\G$ given by $T_{(A,\Psi)}:=\ker(d^{0,*}_{A,\Psi})$. Let $\pi : T_{(A,\Psi)} \to \bonf_{\s}$ be the projection map given by $(a,\Phi) \mapsto [A+a,\Psi + \Phi]$. The induced map
\begin{equation*}
 T_{(A,\Psi)} / \Gamma_{(A,\Psi)} \to \bonf_{\s} \ ,
\end{equation*}
yields a homeomorphism of a neighbourhood of $[0] \in T_{(A,\Psi)}/\Gamma_{(A,\Psi)}$ onto a neighbourhood of $[A,\Psi]$ in the configuration space $\bonf_\s$ of $PU(2)$ monopoles.
\end{proposition}

As we have seen before, a point $[A,0] \in \bonf_\s$ is always a fixed point of the circle action on $\bonf_{\s}$. On the slice $T_{(A,0)}$ we have a circle action $S^1 \times T_{(A,0)} \to T_{(A,0)}$ given by $(z,(a,\Phi)) \mapsto (a, z \Phi)$. The actions of $\Gamma_{(A,0)}$ and $S^1$ commute and they factor through an obvious action of the group $\Gamma_{(A,0)} \times_{\Z/2} S^1$. To simplify the notation we will now write $\Gamma$ instead of $\Gamma_{(A,0)}$. The projection $\pi: T_{(A,0)} \to \bonf_\s$ is then equivariant with respect to the action of $\stab$ on the slice and the $S^1$-action $\rho$ on $\bonf_\s$ of section \ref{circle action}. Therefore we get
\begin{proposition}\label{proposition 2.9}
	The map 
	\begin{equation*}
		T_{(A,0)} / \stab \to \bonf_\s/S^1	 	
	\end{equation*}
	induced from the projection $\pi$ yields a homeomorphism of a neighbourhood of $[0]$ in $T_{(A,0)} / \stab$ onto a neighbourhood of $[A,0]$ in the $S^1$-quotient $\bonf_\s / S^1$.  
\end{proposition}

We will now use this result to describe a neighbourhood of a point $[A,\Psi]$ in the moduli space $M_{\kappa,\s}^{w}$. Let $\mathscr{F}_{\bfomega,\bfalpha}$ be the map given by the left-hand side of the $PU(2)$ monopole equations (\ref{perturbed monopoles}) where the perturbations $(\bfomega,\bfalpha)$ are kept fixed. We consider the restriction of this map to the slice $T_{(A,\Psi)}$:
\begin{equation} \label{map on slice}
	\begin{split}
		f : \ & T_{(A,\Psi)} \to L^2_{l-1}(X;S^- \tensor E \oplus \Lambda^2_+ \tensor \, \su(E)) \\
			& (a,\Phi) \mapsto \mathscr{F}_{\bfomega,\bfalpha}(A + a, \Psi + \Phi)
	\end{split}
\end{equation}
This map is equivariant with respect to the natural action of $\Gamma_{(A,\Psi)}$ on both spaces. Furthermore, it is equivariant with respect to the natural $\stab$ actions in case we consider a fixed point of the circle action of the form $[A,0]$. 
From the above Propositions we therefore get:

\begin{proposition}\label{proposition local model}
	The projection $\pi$ of Proposition \ref{slices} induces a homeomorphism of a neighbourhood of the origin in $f^{-1}(0) / \Gamma_{A,\Psi}$ onto a neighbourhood of $[A,\Psi]$ in $M^{w}_{\kappa,\s}$. Furthermore, it induces a homeomorphism of the origin in $f^{-1}(0) / \stab$ onto a neighbourhood of $[A,0]$ in $M^{w}_{\kappa,\s}/S^1$.
\end{proposition}

Next it can be shown that this description of the neighbourhood of $[A,\Psi]$ in the moduli space by the zero-set of the map $f$ modulo stabiliser can be cut down to the zero-set of a map $h: H^1_{\bfomega,\bfalpha} \to H^2_{\bfomega,\bfalpha}$ between the finite-dimensional cohomology spaces of the elliptic complex (\ref{elliptic complex}) associated to $(A,\Psi)$, modulo stabiliser. See the discussion in \cite[section 4.2.4 and 4.2.5]{DK}. We will now make this local description explicit around fixed points $[A,0]$ of the circle action on the moduli-space. 

Let $Q$ be the derivative of $f$ at $0$. It equals the restriction of $d^1_{\bfomega,\bfalpha}$ of the elliptic complex (\ref{elliptic complex}) to the slice $T_{(A,0)}$. This map $Q$ is a $\stab$-equivariant Fredholm map. Note that $\stab$ acts isometrically on both $T_{(A,0)}$ and the image space of $Q$. Now there is a topological decomposition of the slice as
\begin{equation*}
 	T_{(A,0)} = H^1_{(\bfomega,\bfalpha)} \oplus T' \ ,
\end{equation*}
where $H^1_{(\bfomega,\bfalpha)}$ is the kernel of $Q$ and where $T'$ is a $\stab$-invariant complement of it. Simply take $T'$ to be the $L^2_{l}$-orthogonal complement of $H^1_{(\bfomega,\bfalpha)}$ inside the slice $T_{(A,0)}$. There is also a topological decomposition of the target space as
\begin{equation*}
 	L^2_{l-1}(X;S^- \tensor E \oplus \Lambda^2_+ \tensor \, \su(E))	
			= H^2_{(\bfomega,\bfalpha)} \oplus \im(Q) \ ,
\end{equation*}
where the harmonic space $H^2_{(\bfomega,\bfalpha)}$ is given by  $\ker((d_{A,\bfomega}^+)^*) \oplus \ker(\dirac_{A,\bfalpha}^-)$, which is equally a $\stab$-invariant subspace. Therefore, $Q$ is a map
\begin{equation*}
	 Q: H^1_{(\bfomega,\bfalpha)} \oplus T' \to H^2_{(\bfomega,\bfalpha)} \oplus \im(Q) \ ,
\end{equation*}
and the restriction $Q' : T' \to \im(Q)$ is an equivariant isomorphism of Hilbert spaces. \\

Let $p: H^2_{(\bfomega,\bfalpha)} \oplus \im(Q) \to \im(Q)$ be the orthogonal projection.
\begin{proposition}
There is a $\stab$-equivariant diffeomorphism $g$ of a neighbourhood of the origin in the slice $T_{(A,0)}$ so that 
\begin{equation*} 
 	p \circ f \circ g = p \circ Q \ .
\end{equation*}
\end{proposition}
{\em Proof:} To simplify the notation we shall only write $H^1$ for the harmonic space $H^1_{\bfomega,\bfalpha}$. We define the map 
\begin{equation*}
	G  : H^1 \oplus T' \to H^1 \oplus T' \\
\end{equation*}
by the formula $ G(h,t):= (h,Q'^{-1} \circ p \circ f (h,t))$, where $h \in H^1, t \in T'$. This is a $\stab$-equivariant map. Its derivative at $(0,0)$ is easily seen to be the identity. Therefore $G$ is a diffeomorphism of a neighbourhood of $(0,0)$ onto a neighbourhood of $(0,0)$. Let $g$ be its inverse, which is necessarily equivariant. 
However, the equation $G \circ g (h,t) = (h,t)$ is simply equivalent to 
\begin{equation*}
 p \circ f \circ g (h,t) =  Q'(t) = p \circ Q (h,t) \\ ,
\end{equation*}
which is the equality we sought to prove. \qed

\begin{corollary} \label{finite dimensional model}
Suppose $g$ is a $\stab$-equivariant diffeomorphism as in the last proposition. Then we have
\begin{equation*}
	f \circ g (h,t) = (\alpha(h,t),Q'(t)) \ ,
\end{equation*}
where $\alpha$ is a $\stab$-equivariant map which has vanishing derivative at $0$. As a consequence, up to a $\stab$-equivariant diffeomorphism, the zero-set $f^{-1}(0)$ is given by the zero-set of the $\stab$-equivariant map $\alpha(- , 0) : H^1 \to H^2$. 
\end{corollary}

The following statement is a simple consequence of the implicit function theorem and makes the charts of the moduli space $M^{w,**}_{\s,\kappa}$ a little more explicit:

\begin{lemma} \label{charts}
 Suppose $[A,\Psi]$ is a regular monopole. Then there is a smooth map $h : U \to T_{(A,\Psi)}$, defined on a neighbourhood $U$ of $ 0 \in H^1_{(A,\Psi)}$ which yields a parametrisation of the zero-locus $f^{-1}(0) \subseteq T_{(A,\Psi)}$ of the map (\ref{map on slice}) around $0$. Furthermore, $h$ is of the form 
\begin{equation*}
 h(a,\Phi) = (a,\Phi) + q(a,\Phi) \ ,
\end{equation*}
where the derivative of the map $q$ vanishes at the origin. If $\Psi = 0$ the map $h$ is $S^1$-equivariant.

If $(A,\Psi)$ has trivial stabiliser the composition $\pi \circ h_{(A,\Psi)}$ yields a smooth parametrisation of the moduli space $M^{w,**}_{\s,\kappa}$ in a neighbourhood of $[A,\Psi]$, where $h_{(A,\Psi)}$ denotes the map $(a,\Phi) \mapsto (A,\Psi) + h(a,\Phi)$. 
\end{lemma}



\section{The main Theorem and its proof}
Our main result is the vanishing of the Casson-type invariant defined in Section \ref{definition casson}, Definitions \ref{casson invariant perturbation} and \ref{casson invariant without perturbation}.
\begin{theorem}\label{vanishing theorem}
Suppose the negative definite four-manifold $X$ has first Betti-number $b_1(X)=1$ and admits a class $w \in H^2(X;\Z)$ so that $\langle w^2,[X]\rangle$ is divisible by four, and so that the image of $w$ in $H^2(X;\Z) / \text{\em Torsion}$ is not divisible by 2. 

Under this definition we have defined the Casson-type invariant $n_{w,o}(X)$ in Definitions \ref{casson invariant perturbation} and \ref{casson invariant without perturbation} as the signed count of zero-dimensional moduli space $M_0^w(\bfomega)$ associated to a generic perturbation $\bfomega \in \V$ which was assumed small enough so that the compactness of the main stratum is preserved. Then the equation
\begin{equation*}
 n_{w,o}(X) = 0 \ 
\end{equation*}
always holds.
 \end{theorem}
\subsection{Sketch of proof} Before proving the theorem in the following section we shall first give a sketch of it, omitting the discussion of regularity and orientations. 
\\

Suppose the Casson-invariant is defined by the count $n_{w,o}(X) = \# M_{0}^{w}$.
We will prove the theorem by constructing a suitable compact moduli space $M_{0,\s}^{w}$ of $PU(2)$ monopoles of expected dimension 2, containing the instanton moduli space $M_{0}^{w}$ as the subspace of monopoles with vanishing spinor.

The circle action of section \ref{circle action} restricts to a circle action on the moduli space $M_{0,\s}^{w}$. The instanton moduli space $M_{0}^{w}$ is precisely equal to the fixed-point set of this circle action -- there are no fixed points consisting of abelian $U(1)$ Seiberg-Witten monopoles because the bundle $E \to X$ does not split topologically according to Lemma \ref{no reductions} whereas this would be necessary by the explicit descriptions of these fixed points by Proposition 
\ref{fixed points}. The same is even true for the lower strata $M_{0,\s}^{w}$ which may a priori have non-empty subsets from the Uhlenbeck compactification of $M_{0,\s}^{w}$.

The circle action is free on the complement which is equal to $M_{0,\s}^{w,**}$. This subspace admits the structure of a smooth 2-dimensional manifold for generic perturbation.

As a consequence of our discussion on $S^1$-equivariant local models of the moduli space around fixed points of the $S^1$ action, 
the quotient $M_{0,\s}^{w}/S^1$ of the monopole moduli space by the circle action is a smooth {\em compact} 1-dimensional manifold with boundary. The boundary can be identified with the instanton moduli space $M_{0}^{w}$. Therefore the instanton moduli space $M_{0}^{w}$ is smoothly cobordant to the empty space, and so $n_{w,o}(X) \equiv 0 \text{ mod } 2$. 

As a final step we will discuss orientations. A choice of trivialisation of the determinant line bundle of the deformation operator $-d^*_{A} \oplus d_A^+$ determines an orientation of the instanton moduli space $M^{w}_{0}$. It also determines a natural orientation of the $PU(2)$ monopole moduli space and an orientation of the quotient $M_{0,\s}^{w,**}/S^1$, after introducing a suitable convention for the orientation of the quotient. We will then show that the two orientations of the instanton moduli space $M^{w}_{0}$ -- one as above, and the other as the boundary of the oriented 1-dimensional manifold $M^{w}_{0,\s}/S^1$ -- do in fact coincide. Therefore $M_{0}^{w}$ is smoothly orientedly cobordant to the empty space, and so $n_{w,o}(X) = 0$. 

\subsection{Regularity and the one-dimensional cobordism}
We will choose a $Spin^c$ structure $\s$ according to Proposition \ref{index}. The complex index of the Dirac-operator $D_{A}^+$ is then equal to one. By the discussion in section \ref{properties of monopoles} and the expected dimension $0$ of the instanton moduli space $M^{w}_0$ in Proposition \ref{our moduli space}, the expected dimension of the $PU(2)$ monopole moduli space $M^{w}_{0,\s}$ is two.

Next, we shall make a convenient genericity assumption. Note that countable intersections of residual sets are residual. By Theorem \ref{transversality instantons}, Corollary \ref{generic regularity monopoles}, Lemma \ref{surjectivity dirac} and Theorem \ref{regularity lower strata} there is therefore a residual set of perturbation parameters $(\bfomega,\bfalpha) \in \V \times \W$ and of parameters $(\beta,K,\gamma) \in \mathscr{P}^{l}$ so that the following condition holds:

\begin{condition} \label{genericity assumption}
1. The moduli space of instantons $M_{0}^{w}(\bfomega,h,m)$ with cut off perturbations as in Section \ref{moduli spaces with cut off perturbations} consists of regular points only. It is therefore a compact zero-dimensional manifold and consists of a finite number of points. \\
2. The moduli space $M_{0,\s}^{w}(\bfomega,\bfalpha)$ is compact and the subspace $M_{0,\s}^{w,**}(\bfomega,\bfalpha)$ is a smooth two-dimensional manifold with a free circle-action. \\
3. For each of the finite number of instantons $[A]$ occuring in $M_{0}^{w}(\bfomega,h,m)$ the perturbed Dirac-operator $D^+_{A,\bfalpha}$ has trivial cokernel. 
\end{condition}
In fact, the monopole moduli space $M_{0,\s}^{w}(\bfomega,\bfalpha)$ is compact for the following reason. The lower strata of its Uhlenbeck compactification comprise, according to Section \ref{compactification of monopoles}, elements in the monopole moduli spaces $M_{-n,\s}^{w}$ for $n \geq 1$ (with zero holonomy perturbation). However, these moduli spaces do not contain fixed points of the circle action by Lemma \ref{no reductions}, Proposition \ref{fixed points}, and because the instanton moduli spaces $M^{w}_{-n}$ are empty for $n \geq 1$. Therefore $M_{-n,\s}^{w} = M_{-n,\s}^{w,**}$. But these moduli spaces have strictly negative expected dimension by Proposition \ref{expected dimensions compactification}, so they must be generically empty by Theorem \ref{regularity lower strata} with regard to the `auxiliary perturbation parameters' $(\beta,K,\gamma)$. 
\\

The Casson-invariant is the signed count $n_{w,o}(X):= \# M_{0}^{w}(\bfomega) = \# M_{0}^{w}(\bfomega,h,m)$ according to Proposition \ref{well-definedness} and the remarks in Section \ref{moduli spaces with cut off perturbations}. In order to keep the exposition easier we shall now simply write $M_{0}^{w}(\bfomega)$ when we actually mean the moduli space $M_{0}^{w}(\bfomega,h,m)$ with cut off perturbations. \\ 

Note that the circle action $ \rho $ of section \ref{circle action} restricts to a circle action on the monopole moduli space $M_{0,\s}^{w}(\bfomega,\bfalpha)$. According to Proposition \ref{fixed points} and the it following remark, the fixed point set of this restriction equals the instanton moduli space,
\begin{equation*}
 	M_{0,\s}^{w}(\bfomega,\bfalpha)^{S^1} \cong M_{0}^{w}(\bfomega) \ ,
\end{equation*}
where an instanton $[A]$ is considered as $PU(2)$ monopole $[A,0]$. Note also that the complement of the fixed-point set of the circle action consists of monopoles with irreducible connection and non-vanishing spinor,
\begin{equation*}
 	M_{0,\s}^{w}(\bfomega,\bfalpha) \setminus M_{0}^{w}(\bfomega) = M_{0,\s}^{w,**}(\bfomega,\bfalpha) \ .
\end{equation*}
This complement $M_{0,\s}^{w,**}(\bfomega,\bfalpha)$ is a 2-dimensional smooth $S^1$-space, and it has a natural compactification (inside the entire moduli space $M_{0,\s}^{w}(\bfomega,\bfalpha)$) with ends given by the instanton moduli space $M_{0}^{w}(\bfomega)$. However, we want to be sure that each instanton only corresponds to one `end' of the moduli space $M_{0,\s}^{w,**}(\bfomega,\bfalpha)$. For this we have to study the local structure of the moduli space $M_{0,\s}^{w}(\bfomega,\bfalpha)$ in the neighbourhood of an instanton $[A,0]$. \\

By Proposition \ref{proposition local model} and Corollary \ref{finite dimensional model} a neighbourhood of an instanton $[A,0]$ in the monopole moduli space $M_{0,\s}^{w}(\bfomega,\bfalpha)$ or its $S^1$-quotient is described be a $\stab$-equivariant map
\begin{equation*}
	\alpha(-,0) : H^1_{(A,0),(\bfomega,\bfalpha)} \to H^2_{(A,0),(\bfomega,\bfalpha)} \ ,
\end{equation*}
where the spaces $H^i_{(A,0),(\bfomega,\bfalpha)}$, $i=1,2$ are the cohomology spaces (or harmonic spaces) of the elliptic deformation complex (\ref{elliptic complex}) associated to the solution $(A,0)$ of the perturbed monopole equations (\ref{perturbed monopoles}), and are given by
\begin{equation*}
 	\begin{split}
 	 	H^1_{(A,0),(\bfomega,\bfalpha)} & = \ker(\delta_{A,\bfomega}) \oplus 	\ker(\dirac_{A,\bfalpha}^+) \\
		H^2_{(A,0),(\bfomega,\bfalpha)} & = \coker(d_{A,\bfomega}^+) \oplus \coker(\dirac^+_{A,\bfalpha}) \ .
 	\end{split}
\end{equation*}
Under the above genericity assumption these spaces are given by
\begin{equation*}
	\begin{split}
	H^1_{(A,0)} & = \ker(\dirac^+_{A,\bfalpha}) \cong \C \\
	H^2_{(A,0)} & = 0 \ ,
	\end{split}
\end{equation*}
where we have omitted the dependence on the perturbation.

The $S^1$-action on $H^1_{(A,0)} \subseteq T_{(A,0)}$ corresponding via Proposition \ref{proposition 2.9} to the action $\rho$ on the configuration space $\bonf^{w}_{0,\s}$ is so that the $S^1$-equivariant identification $\ker(\dirac^+_{A,\bfalpha}) \cong \C$ corresponds to the standard action of $S^1$ on $\C$. 

Note that the stabiliser $\Gamma_{(A,0)}$ is isomorphic to $\Z/2$. Proposition \ref{proposition local model} now implies that a neighbourhood of $[A,0]$ in the moduli space $M_{0,\s}^{w}(\bfomega,\bfalpha)$ is homeomorphic to a neighbourhood of $0$ in the 
quotient $\C/(\Z/2)$ - it can be taken as a cone over the real projective space $\mathbb{RP}^1$. It also implies that a neighbourhood of $[A,0]$ in the $S^1$-quotient $M_{0,\s}^{w}(\bfomega,\bfalpha)/S^1$ is homeomorphic to a neighbourhood of $0$ in the quotient $\C/S^1 \cong [0,\infty) \subseteq \R$. Thus we get

\begin{proposition}\label{structure of quotient}
 The $S^1$-quotient of the monopole moduli space $M_{0,\s}^{w}(\bfomega,\bfalpha)$ is a smooth one-dimensional manifold with boundary. Its boundary can be identified with the instanton moduli space $M_{0}^{w}(\bfomega)$. 
\end{proposition}
\begin{corollary}
The instanton moduli space $M_{0}^{w}(\bfomega)$ is cobordant to the empty space. As a consequence, we must have $n_{w,o}(X) \equiv 0 \text{ mod } 2$.
\end{corollary}

\subsection{Consideration of orientations}
We shall first agree that the boundary of an oriented manifold is oriented by the `outward normal first' convention. In this way, for instance, the one-sphere $S^1$, seen as the boundary of the unit disc in $\C$ with its complex orientation, has its orientation `counterclock-wise'. 

We shall continue to suppose that perturbations $(\bfomega,\bfalpha) \in \V \times \W$ are chosen so that the Condition \ref{genericity assumption} holds. However, we shall not make this explicit in the notation of the harmonic spaces $H^1_{(A,\Psi),(\bfomega,\bfalpha)} = \ker(D_{(A,\Psi),(\bfomega,\bfalpha)})$ anymore and simply write $H^1_{(A,\Psi)}$ instead, and likewise for the deformation operators.
\\

The zero-dimensional moduli space $M_0^w(\bfomega)$ is oriented by a choice of trivialisation $o$ of the determinant line bundle $\Lambda_0 \to \bonf^{w}_0$, see section \ref{orientations instantons}. 

According to Corollary \ref{natural orientation monopoles}, the regular subspace $M^{w,**}_{\s,0}(\bfomega,\bfalpha)$ of the $PU(2)$ monopole moduli space is naturally oriented by a choice of trivialisation of the same determinant line bundle $\Lambda_0 \to \bonf^{w}_0$, and we shall choose the same trivialisation $o$ as above. 

The circle $S^1$ acts freely and smoothly on $M^{w,**}_{\s,0}(\bfomega,\bfalpha)$.
An orientation of the $S^1$ quotient $M^{w,**}_{\s,0}(\bfomega,\bfalpha)/S^1$ is fixed by the following convention. Let $[A,\Psi]$ belong to $M^{w,**}_{\s,0}(\bfomega,\bfalpha)$, and let $[[A,\Psi]]$ be the corresponding element in the $S^1$ quotient. We require that the direct sum of tangent spaces
\begin{equation} \label{orientation convention}
	T_{[A,\Psi]}M^{**} \, = \, T_{1} S^1.[A,\Psi] \, \oplus \, T_{[[A,\Psi]]}M^{**}/S^1 
\end{equation}
is a direct sum of oriented vector spaces, where $S^1.[A,\Psi]$ denots the $S^1$ orbit through $[A,\Psi]$. In other words, an oriented basis of 
$T_{[A,\Psi]}M^{**}$ is obtained by completing an oriented basis of the orbit with an oriented basis of $T_{[[A,\Psi]]}M^{**}/S^1$.
\\

On the other hand, according to Proposition \ref{structure of quotient}, the instanton moduli space $M^{w}_0(\bfomega)$ 
can be seen as the boundary of the 1-dimensional manifold $M^{w}_{\s,0}(\bfomega,\bfalpha)/S^1$. The open submanifold $M^{w,**}_{\s,0}(\bfomega,\bfalpha)/S^1$ was given an orientation in the last paragraph.  This induces an orientation on the boundary $M^{w}_0(\bfomega)$. The proof of Theorem \ref{vanishing theorem} will be complete if we show that the two orientations on $M^{w}_0(\bfomega)$ coincide, because then the oriented moduli space $M^w_0(\bfomega)$ is orientedly cobordant to the empty space and therefore the Casson-type invariant $n_{w,o}(X)$ must be zero. 
\\

The harmonic space $H^1_{(A,0)} = \ker(\dirac_{A,\bfalpha}^+)$ is oriented by the trivialisation $o$ of $\Lambda_0$, again by Corollary \ref{natural orientation monopoles}. If the instanton $[A] \in M^{w}_0(\bfomega)$ has orientation $+1$ the space $H^1_{(A,0)}$ is oriented by its natural complex orientation, otherwise it is oriented with the opposite of the complex orientation. 
\\

Now let us suppose that the monopole $[A',\Psi']$ is regular and has trivial stabiliser, and that further $(A',\Psi') \in T_{(A,0)}$ is close enough to $(A,0)$ so that it is in the range of the parametrisation $h_{(A,0)}$ of Lemma \ref{charts}, and let $h_{(A,0)}(\Phi) = (A',\Psi')$ for $\Phi \in H^1_{(A,0)}$. The tangent space $T_{\Phi}H^1_{(A,0)}$ of $H^1_{(A,0)}$ at $\Phi$ is canonically identified with $H^1_{(A,0)}$. Let us denote explicitly the following two cases:

\begin{enumerate}
\item If $H^1_{(A,0)}$ admits its complex orientation the tangent space $T_{\Phi}H^1_{(A,0)}$ has an oriented basis consisting of
\begin{equation*}
	\left(\frac{\Phi}{\norm{\Phi}}, \frac{i \, \Phi}{\norm{\Phi}}\right) \ .
\end{equation*}
\item If $H^1_{(A,0)}$ admits the opposite of the complex orientation the tangent space $T_{\Phi}H^1_{(A,0)}$ has an oriented basis consisting of
\begin{equation*}
	\left(-\frac{\Phi}{\norm{\Phi}}, \frac{i \, \Phi}{\norm{\Phi}}\right) \ .
\end{equation*}
\end{enumerate}

\begin{lemma} \label{orientation preserving}
	Let $h: U \to T_{(A,0)}$ be a parametrisation of the zero-locus of the monopole map restricted to the slice (\ref{map on slice}) as in Lemma \ref{charts}. After possibly restricting $h$ to a smaller neighbourhood $U'$ of $0 \in H^1_{(A,0)}$,
	the composition 
	\begin{equation*}
	\pi \circ h_{(A,0)}|_{U' \setminus \{ 0 \}} : U' \setminus \{ 0 \} \to M^{w,**}_{\s,0}
	\end{equation*}
	is a local diffeomorphism that is two-to-one and that is orientation-preserving.
%
\end{lemma}
Although this lemma seems obvious our proof is slightly technical. We defer it to the end of this section. \\

We will now complete the proof of Theorem \ref{vanishing theorem} assuming Lemma \ref{orientation preserving}. Without loss of generality we may assume that the neighbourhood $U$ equals the neighbourhood $U'$ of the last lemma. 

As $h$ is $S^1$-equivariant the vector 
\begin{equation*}
 	d_{\Phi} (\pi \circ h_{(A,0)}) \left(\frac{i \, \Phi}{\norm{\Phi}}\right) \, \in \, T_{[A',\Psi']}M^{w,**}_{\s,0}
\end{equation*}
yields a positive basis of the tangent space to the $S^1$-orbit $S^1 . [A',\Psi']$ at $[A',\Psi']$. 
\\

Suppose the instanton $[A] \in M^{w}_{0}(\bfomega)$ counts as $+1$. By Lemma \ref{orientation preserving} the two vectors 
\begin{equation*}
 \left(d_{\Phi} (\pi \circ h_{(A,0)})  \left(\frac{\Phi}{\norm{\Phi}}\right), d_{\Phi} (\pi \circ h_{(A,0)}) \left(\frac{i \, \Phi}{\norm{\Phi}}\right)\right)
\end{equation*}
form a positive basis of the tangent space $T_{[A',\Psi']}M^{**}$. By our orientation convention (\ref{orientation convention}) a positive basis of the tangent space
$T_{[[A',\Psi']]}M^{**}/S^1$ to the quotient $M^{**}/S^1$ is given by the vector 
\begin{equation*}
 - d_{\Phi} (\pi \circ h_{(A,0)}) \left(\frac{ \Phi}{\norm{\Phi}}\right) \ .
\end{equation*}
This vector points `outwards' towards the boundary instanton $[A]$. Therefore the orientation of $[A]$ as boundary point of the one-dimensional manifold with boundary $M^{w}_{\s,0}/S^1$ is positive. 

Now suppose the instanton $[A] \in M^{w}_{0}(\bfomega)$ counts as $-1$. Now the two vectors \begin{equation*}
 \left(-d_{\Phi} (\pi \circ h_{(A,0)})  \left(\frac{\Phi}{\norm{\Phi}}\right), d_{\Phi} (\pi \circ h_{(A,0)}) \left(\frac{i \, \Phi}{\norm{\Phi}}\right)\right)
\end{equation*}
form a positive basis of the tangent space $T_{[A',\Psi']}M^{**}$ and a positive basis of the tangent space $T_{[[A',\Psi']]}M^{**}/S^1$ to the quotient $M^{**}/S^1$ is given by the vector 
\begin{equation*}
  d_{\Phi} (\pi \circ h_{(A,0)}) \left(\frac{ \Phi}{\norm{\Phi}}\right) \ .
\end{equation*}
This vector points `inwards' away from the boundary instanton $[A]$. Therefore the orientation of $[A]$ as boundary point of the one-dimensional manifold with boundary $M^{w}_{\s,0}/S^1$ is negative. \qed 

\subsection{Proof of Lemma \ref{orientation preserving}}
By previous results we will only have to check that the map $\pi \circ h_{(A,0)}|_{U \setminus \{ 0 \}}$ indeed has an orientation-preserving derivative for sufficiently small neighbourhoods $U$ of $0 \in H^1_{(A,0)}$. 
\\

Suppose $h_{(A,0)}(\Phi) = (A',\Psi') \in (A,0) + T_{(A,0)}$. Let $h_{(A',\Psi')} : H^1_{(A',\Psi')} \to (A',\Psi') + T_{(A',\Psi')}$ be a parametrisation of $M^{w,**}_{\s,0}$ around $[A',\Psi']$ as in Lemma \ref{charts}. In order to prove the lemma, we will need a convenient restriction of the map $\pi \, \circ \, h_{(A,0)}$ to a neighbourhood $V$ of $\Phi \in H^1_{(A,0)}$ through the chart $(\pi \, \circ \, h_{(A',\Psi')})^{-1}$. This map 
\begin{equation} \label{kartenwechsel}
	k:=(\pi \circ h_{(A',\Psi')})^{-1} \circ (\pi \circ h_{(A,0)}|_{V}) \, : \, V \subseteq H^1_{(A,0)} \, \to \, H^1_{(A',\Psi')}
\end{equation}
can be described in an alternative way. 
\\ 

The spaces $H^1_{(A,0)}$ and $H^1_{(A',\Psi')}$ lie in the different slices $T_{(A,0)}$ and $T_{(A',\Psi')}$ to the action of the gauge group $\G$ on $\conf$. We wish to define a `gauge fixing changing' map 
\begin{equation*}
 g: (A,0) + T_{(A,0)} \to (A',\Psi') + T_{(A',\Psi')} \ ,
\end{equation*}
at least in a neighbourhood of $(A',\Psi')$, with the following two properties: 
\begin{enumerate}
 \item $g((A',\Psi')) = (A',\Psi')$ \ ,
 \item $[g((A',\Psi') + (a,\Sigma))] = [(A',\Psi') + (a,\Sigma)] $ .
\end{enumerate}
\begin{lemma}\label{eichfixierung}
 There exists a smooth gauge fixing changing map $g$ satisfying the above two properties so that its derivative $d_{(A',\Psi')} g : T_{(A,0)} \to T_{(A',\Psi')}$ at $(A',\Psi')$ has the following shape
\begin{equation*}
	(d_{(A',\Psi')} g) \, (a,\Sigma) = (a,\Sigma) + c(a,\Sigma) \ ,
\end{equation*}
where the norm of the linear map $c$ can be made as small as we wish by choosing $(A',\Psi')$ close enough to $(A,0)$. 
\end{lemma}
{\em Proof:}
We apply the implicit function theorem to the map
\begin{equation*}
 \begin{split}
  	F: L^2_{l+1}(X;\su(E)) \times T_{(A,0)} & \to L^2_{l-1}(X;\su(E)) \\
	(\zeta,a,\Sigma) & \mapsto d_{(A',\Psi')}^{0,*} \, (\exp(\zeta)(A'+a,\Psi'+\Sigma)) \ , 
 \end{split}
\end{equation*}
where $\exp: L^2_{l+1}(X;\su(E)) \to \G$ denotes the exponential from the Lia-algebra of the gauge group to the gauge group. Note the following two partial derivatives of $F$ at $(0,(0,0))$:
\begin{equation}\label{partielle ableitungen}
\begin{split}
 \left. \frac{\partial F}{\partial \zeta} \right|_{(0,(0,0))} (\zeta) & = d_{(A',\Psi')}^{0,*} d_{(A',\Psi')}^{0} (\zeta) = \Delta^0_{(A',\Psi')} (\zeta) \\
\left. \frac{\partial F}{\partial (a,\Sigma)} \right|_{(0,(0,0))} (a,\Sigma) & = 
	d_{(A',\Psi')}^{0,*} (a,\Sigma) \ . 
\end{split}
\end{equation}
As $(A',\Psi')$ is irreducible by assumption, the Laplacian $\Delta^0_{(A',\Psi')} = d_{(A',\Psi')}^{0,*} d_{(A',\Psi')}^{0}$ is an isomorphism. Therefore, by the implicit function theorem there is a map $\gamma$ from a neighbourhood of $0$ in $T_{(A,0)}$ to a neighbourhood of $0$ in $L^2_{l+1}(X;\su(E))$ so that 
\begin{equation*}
  F(\gamma(a,\Sigma),(a,\Sigma)) = 0 \  \ 
\end{equation*}
for alle $(a,\Sigma)$ in the domain of definition of $\gamma$. Equivalently,
\begin{equation*}
  d_{(A',\Psi')}^{0,*}(\exp(\gamma(a,\Sigma)) (A'+a,\Psi'+\Sigma)) = 0 \ 
\end{equation*}
for all these $(a,\Sigma)$. 
The required map $g$ is then given by 
\[g((A',\Psi') + (a,\Sigma)) := \exp(\gamma(a,\Sigma))(A'+a,\Psi'+\Sigma) \ . 
\]
For the derivative we have
\begin{equation*}
 (d_{(A',\Psi')} g) \, (a,\Sigma) = 
	(a,\Sigma) - \, d^0_{(A',\Psi')} \left(\Delta^0_{(A',\Psi')}\right)^{-1} d_{(A',\Psi')}^{0,*} \, (a,\Sigma) \ .
\end{equation*}
Indeed, it is immediate to see that 
\begin{equation*}
 (d_{(A',\Psi')} g) \, (a,\Sigma) = 
	(a,\Sigma)  +  d^0_{(A',\Psi')} \circ \left. \frac{\partial \gamma}{\partial (a,\Sigma)} \right|_{(0,0)} (a,\Sigma) \ .
\end{equation*}
On the other hand, by differentiating the map 
\[
	(a,\Sigma) \mapsto F(\gamma(a,\Sigma),(a,\Sigma)) \ ,
\]
which is constant $0$, one sees that 
\[
 0 =  \left. \frac{\partial F}{\partial \zeta} \right|_{(0,(0,0))} \circ \left. \frac{\partial \gamma}{\partial (a,\Sigma)} \right|_{(0,0)} + 
\left. \frac{\partial F}{\partial (a,\Sigma)} \right|_{(0,(0,0))} \ . 
\]
The claimed formula for $(d_{(A',\Psi')} g) \, (a,\Sigma)$ now follows from 
the two formulae (\ref{partielle ableitungen}) above, and the fact that the first of these partial derivatives of $F$, given by the Laplacian $\Delta^0_{(A',\Psi')}$, is an isomorphism.

Now notice that the operators $d^0_{(A',\Psi')}$ and $\Delta^0_{(A',\Psi')}$ vary continously with $(A',\Psi')$. But for $(a,\Sigma) \in T_{(A,0)} = \ker (d_{(A,0)}^{0,*})$ we have 
\begin{equation*}
 d_{(A',\Psi')}^{0,*} \, (a,\Sigma) = (d_{(A',\Psi')}^{0,*} \, - d_{(A,0)}^{0,*}) \, (a,\Sigma) \ ,
\end{equation*}
and the operator $d_{(A',\Psi')}^{0,*} \, - d_{(A,0)}^{0,*}$ can be made as small as we wish by chosing $(A',\Psi')$ close enough to $(A,0)$. \qed

Now notice that the map $k$ of (\ref{kartenwechsel}), at least when restricted to a sufficiently small neighbourhood of $\Phi \in H^1_{(A,0)}$, is given by
\begin{equation*}
 	k (0,\Phi+ \Sigma) = h_{(A',\Psi')}^{-1}((g \circ h_{(A,0)})(0,\Sigma)) \ .
\end{equation*}
By Lemmata \ref{eichfixierung} and \ref{charts} the derivative of $k$ at $\Phi$ has the form 
\begin{equation*}
 d_{\Phi} k \, (\Sigma) = \Sigma + c_{(A',\Psi')} (\Sigma) \ ,  
\end{equation*}
where the linear map $c_{(A',\Psi')}$ can be made as small as we like by chosing $(A',\Psi')$ close enough to $(A,0)$. In particular, for $\norm{c_{(A',\Psi')}} < 1$ the map $d_{\Phi} k$ is an isomorphism. 
\\

We will compare the derivative $d_{\Phi} k : H^1_{(A,0)} \to H^1_{(A',\Psi')}$ to an orientation preserving map. Let $R$ be a right-inverse to the deformation operator $D_{(A,0)}$. We define linear maps
\begin{equation*}
	\begin{split}
 	P_{(A',\Psi')} : H^1_{(A',\Psi')} & \to H^1_{(A,0)} \\
			(a,\Sigma) & \mapsto (a,\Sigma) + R \,  (D_{(A,0)} - D_{(A',\Psi')}) (a,\Sigma)
	\end{split}
\end{equation*}
which obviously are isomorphisms for $(A',\Psi')$ close enough to $(A,0)$. In a small enough neighbourhood of $(A,0)$ in the configuration space $\conf$, consisting of regular  configurations with at most finite stabilisers only, the maps $P_{(A',\Psi')}$ yield a local trivialisation of the bundle of kernels of the deformation operators $D_{(\bfomega,\bfalpha)}$. In particular, the maps $P_{(A',\Psi')}$ are all orientation-preserving with orientations determined by the (lift of the) determinant line bundle $\Theta_{(\bfomega,\bfalpha)}(1)$ formed by the family of deformation operators. 

We now see that the difference of the two isomorphisms \[ P_{(A',\Psi')}^{-1}\, , d_{\Phi} k \ : \, H^1_{(A,0)} \to H^1_{(A',\Psi')} \] can be made small enough by choosing $(A',\Psi')$ close enough to $(A,0)$. Then both of these maps must be orientation-preserving. \qed

\section{Perspectives}
In this section we shall list some perspectives and potential applications of our main result. Some of these are of a rather speculative nature. 

\subsection{Considerations on four-manifolds with boundary}
The definition of the Casson-type invariant discussed in this article can be extended in various ways to four-manifolds with boundary. Such an extended version of the invariant is more likely not to vanish. For instance, in \cite{LZ} Andrew Lobb and the author have constructed, in a somewhat ad-hoc way, 4-manifolds $X$ with $b_1(X)=1$, $b_2(X)=4$, whose boundary $Y$ are the $0$-surgery of a knot, coming with representations $\rho: \pi_1(X) \to SO(3)$ which restrict to the non-trivial $SO(3)$--representations on four tori $T_1, \dots, T_4$, which span $H_2(X;\Z)$. 

Presumably the easiest way (at least, analytically speaking) to establish an extension of our Casson-type invariant is  to consider negative definite four-manifolds $X$ that have boundary $Y = \partial X$ a collection of integer homology three-spheres. Following ideas of Floer homology \cite{Donaldson_Floer}, one considers a flat connections $\rho$ on the boundary $Y$, so that the expected dimension of the moduli space of projectively flat instantons in a similar bundle $E \to X$ as defined above, with limit (on cylindrical ends) equal to the flat connection $\rho$, is zero. The number $n_{w,o}(X,\rho)$ would then be a count of this moduli space.  

A sample application could now be of the following kind: Assume $n_{w,o}(X,\rho)\neq 0$. Then $Y$ does not bound a positive definite four-manifold $W$ with $b_1(W) = 0$ so that there is only one flat connection $B_\rho$ on the trivial $SU(2)$ bundle on $W$ that extends $\rho$. The proof would use our vanishing result on the closed negative 4-manifold $X \cup_Y \overline{W}$ together with some degeneration arguments. Presumably one would get transversality by holonomy perturbations on $X$ keeping the connections flat on the side $\overline{W}$. This may even work if $\rho$ is the trivial connection on $Y$. A conclusion from $n_{w,o}(X,\rho) \neq 0$ might then be that $Y$ does not bound a simply connected positive definite four-manifold.  

Similar ideas may apply if the boundary manifold $Y$ is not required to be a collection of homology spheres. 

\subsection{Aspherical four-manifolds, four-dimensional Poincar\'e duality groups }
Our vanishing result might possibly be used as an obstruction for a class of four-dimensional Poincar\'e duality ($\text{\em PD-4}$) groups to be the fundamental group of a smooth aspherical four-manifold. 

More precisely, if there were a group $G$ so that $K(G,1)$ admitted the structure of a smooth negative definite 4-manifold $X$ with $b_1(X) = 1$ and $b_2(X) \geq 4$ the invariant $n_{w,o}(X)$ of section \ref{definition casson} associated to a class $w \in H^2(X;\Z)$ so that $\langle w^2,[X]\rangle$ is divisible by four, and so that the image of $w$ in $H^2(X;\Z) / \text{\em Torsion}$ is not divisible by 2, would be zero by our vanishing result. The author has to admit that he has no idea whether such a group $G$ exists at all. 

The {\em unperturbed} moduli space $M^w_0(X)$ consists of equivalence classes of connections $A \in \A(\su(E))$ that are flat. Because $M^w_0(X)$ is the quotient of anti-self-dual connections by the gauge group $\G$, which consists of automorphisms of $E$ with determinant $1$, there is the residual action of $H^1(X;\Z/2)$ on $M^w_0(X)$ with quotient homeomorphic (via the holonomy representation) to the space
\[
	\mathscr{R}_{w}(\pi_1(X)) = \Hom(\pi_1(X),SO(3)) / \SO(3) \ , 
\]
where the action of $SO(3)$ on $\Hom(\pi_1(X),SO(3))$ is by conjugation on the image. Of course, the space $\Hom(\pi_1(X),SO(3)) = \Hom(G,SO(3))$ for a $K(G,1)$ space $X$ is entirely determined by the group $G$. Even more is true: Whether a flat connection $A$ represents a regular point of the moduli space $M^w_0(X)$ is entirely determined by the properties of the group $G$ and the holonomy representation $\rho_A$ determined by $A$. We will outline this shortly in the following interlude.
\\

The deformation complex of a flat connection $A$, seen as instanton representing a point in $M^w_0(X)$, is given by 
\begin{equation}\label{instanton deformation complex}
	0 \to \Omega^0(X;\su(E)) \stackrel{d_A}{ \to} \Omega^1(X;\su(E))
\stackrel{d_{A}^+}{\to} \Omega^2_+(X;\su(E)) \to 0 \ .
\end{equation}
It is regular if the cokernel of $d_A^+$ vanishes. 
For the flat connection $A$ we also have an elliptic complex given by the twisted de Rham complex: 
\begin{equation}\label{twisted de Rham complex}
\begin{split}
	0 & \to \Omega^0(X;\su(E)) \stackrel{d_A}{ \to} \Omega^1(X;\su(E))
\stackrel{d_{A}}{\to} \Omega^2(X;\su(E)) \\ &  \stackrel{d_{A}}{ \to} \Omega^3(X;\su(E))
\stackrel{d_{A}}{ \to} \Omega^4(X;\su(E))
\to 0  \ .
\end{split}
\end{equation}
The first complex (\ref{instanton deformation complex}) gives rise to cohomology groups $H^0_A(X)$, $H^1_A(X)$ and $H^2_A(X)$. Let us denote by $K^0_A(X)$, $K^1_A(X)$, \dots the cohomology groups of the second complex  (\ref{twisted de Rham complex}).
The following proposition is a well-known and easily proved folklore result at least for the untwisted case, and the proof for the twisted case is completely analogous.
\begin{proposition}
Let $A$ be a flat $SO(3)$ connection in the vector bundle $\su(E)$. For the  cohomology spaces in the
above complexes (\ref{instanton deformation complex}) and (\ref{twisted de Rham complex}) we have
\begin{equation*}
	H^0_A(X) \cong K^0_A(X) \ \ \text{ and } \ \ H^1_A(X) \cong K^1_A(X) \ .
\end{equation*}	
\end{proposition}

\begin{remark}
This result generalises to any flat connection in an Euclidean vector bundle which is compatible with the
Euclidean structure. 
\end{remark}

A flat connection $A$ represents a regular point in the moduli space $M^w_0(X)$ if and only if the group $H^1_A(X)$ vanishes. In fact, $H^0_A(X) = 0$ because $A$ is irreducible, and the index of the elliptic complex (\ref{instanton deformation complex}) is zero. Therefore, $H^1_A(X) = 0$ implies $H^2_A(X) = 0$ which we have defined as being `regular' in Section \ref{basic definition of the moduli space} above.

 On the other hand, the cohomology group $H^1_A(X)$, seen as a cohomology group of the twisted de Rham complex, is isomorphic to the first cohomology group $H^1(X;\su(2)_{\ad(\rho_A)})$ of $X$ with local coefficients determined by the action of $\pi_1(X)$ on $\su(2)$ via $\rho$ composed with the adjoint action of $SO(3)$ on its Lie algebra. Whether the latter group vanishes is a matter of the group $G$ and the representation $\rho_A$ alone.
\\

There are analogous topological settings where the action of $H^1(X;\Z/2)$ on the moduli space $M^w_0(X)$ is always free, see for instance \cite[Proposition 4.1 and Section 10.4]{RS}. If we were in such a setting, a sample application of our vanishing result could follow these lines:
Assume we have a $\text{ \em PD-4}$ group $G$ with $b_1(G)=1$ and with diagonalisable negative definite intersection form. Assuming further that  there were a single (or an odd number of) representations $\rho: G \to SO(3)$ with $w_2(\rho)$ equal to the reduction modulo 2 of a class $w \in H^2(G;\Z)$ as above, and that all the twisted cohomology groups $H^1(X;\su(2)_{\ad(\rho)}$ vanished, we would necessarily get a non-zero number $n_{w,o}(X)$ if $G$ were the fundamental group of a smooth aspherical four-manifold $X$. Our vanishing result would therefore contradict this smoothability. 
\\

To get an idea of what phenomena may occur, we shall recall a few facts about Poincar\'e duality groups in higher dimensions: In general, not every Poincar\'e duality group arises as the fundamental group of a closed  aspherical manifold. For instance, there are Poincar\'e duality groups that are not finitely presentable, so cannot arise from a closed manifold. So the question about representability must be asked for finitely presentable groups. In dimension $\geq 5$ techniques from surgery theory may be applied to these problems. There may also be a gap between realisability by a topological/PL/smooth manifold. Davis and Haussmann prove in \cite{DH} that there are topological aspherical manifolds that are not homotopy equivalent to a PL manifold in dimension $\geq 8$, and that there are aspherical PL manifolds not homotopy equivalent to a smooth manifold in dimension $\geq 13$. Our intention is to provide an obstruction to realisability by a smooth manifold in dimension $4$, using gauge theory. Of course, the realisability by a topological manifold is a question in its own right, and the gauge theoretic method cannot be expected to say anything here.

In dimension $4$ Prasad and Yeung's fake projective planes \cite{PY} have to be mentioned. These are smooth aspherical four-manifolds that are rational homology complex projective planes.

\subsection{Teleman's program on class $\text{\em VII}$ surfaces}\label{Teleman's program}
A class $\text{\em VII}$ surface is a closed complex surface with first Betti number $b_1=1$ and with Kodaira dimension $- \infty$. In particular, these are non-algebraic and non K\"ahler surfaces. Surfaces of this class with $b_2=0$ are completely classified, they are either Hopf surfaces or Inoue surfaces.

The current main conjecture in the classification of closed complex surfaces of class $\text{\em VII}$ with positive second Betti number is that any such surface which is {\em minimal} admits a `global spherical shell' -- this is an embedded non-seperating three-sphere admitting a neighbourhood that is biholomorphic to a neighbourhood of the unit sphere in $\C^2$. Such surfaces are well understood by a result of Kato \cite{Kato}. By a result of Dloussky, Oeljeklaus and Toma, a minimal class $\text{\em VII}$ surface with second Betti number $b_2$ admits a global spherical shell if it admits $b_2$ rational curves. Teleman has proved in \cite{T2} that any minimal class $\text{\em VII}$ surface with $b_2=1$ admits a rational curve, using Donaldson theory. 
\\

We will now give a sketch of argument why the non-vanishing of $n_{w,o}(X)$ for a class ${\text{\em VII}}$ surface $X$ with $b_1(X) = 1$ and $b_2(X) \geq 4$ and cohomology class $w \in H^2(X;\Z)$ so that $\langle w^2,[X]\rangle$ is divisible by four, and so that the image of $w$ in $H^2(X;\Z) / \text{\em Torsion}$ is not divisible by 2 would have contradicted the global spherical shell conjecture. In this sense, our vanishing result is a positive result in view of the classification of class $\text{\em VII}$ surfaces. 

If $X$ had a global spherical shell, we would have a connected sum decomposition $X \cong X_0 \, \# \, S^1 \times S^3$. Presumably there would be a perturbation $\bfomega$ so that in the sum $V_{\bfomega}$ of the holonomy perturbation all terms vanish that are associated to loops whose image lies at least in part in the $S^1 \times S^3$ summand, and so that the moduli space $M^w_0(\bfomega)$ is still regular. Any instanton $[A] \in M^w_0(\bfomega)$ would then be flat over the $S^1 \times S^3$ summand and in particular would be trivial over the area where $X_0$ and $S^1 \times S^3$ are glued together. Therefore, any such instanton $[A]$ would extend to an instanton $[A_0]$ over $X_0$, defining an element in a regular moduli space $M^w_0(\bfomega)$ over $X_0$. However, this moduli space has expected dimension equal to $-1$, so it is empty for a perturbation that makes it regular. Therefore, a ${\text{\em VII}}$ surface $X$ with $b_1(X) = 1$ and $b_2(X) \geq 4$ could not have had a non-zero invariant $n_{w,o}(X)$.

\end{document}